\documentclass[12pt]{article}
\usepackage{amsmath}
\usepackage{amssymb}
\usepackage{amsthm}
\usepackage[totalwidth=17cm, totalheight=25cm]{geometry}
\usepackage{graphicx}
\usepackage{mathabx}
\usepackage{tabularx}
	\newcolumntype{C}[1]{>{\centering\arraybackslash}m{#1}} 
	\newcolumntype{R}[1]{>{\raggedleft\arraybackslash}m{#1}} 

\newtheoremstyle{boldplain}
{9pt}
{9pt}
{\itshape}
{}
{\bfseries}
{.}
{.5em}
{\thmname{#1}\thmnumber{ #2}\thmnote{ (#3)}}%

\newtheoremstyle{bolddefinition}
{9pt}
{9pt}
{}
{}
{\bfseries}
{.}
{.5em}
{\thmname{#1}\thmnumber{ #2}\thmnote{ (#3)}}%

\theoremstyle{boldplain}

\newtheorem{cor}[equation]{Corollary}

\newtheorem{lem}[equation]{Lemma}
\newtheorem{lemma}[equation]{Lemma}
\newtheorem{prop}[equation]{Proposition}

\newtheorem{thm}[equation]{Theorem}

\theoremstyle{bolddefinition}
\newtheorem{dfn}[equation]{Definition}
\newtheorem{definition}[equation]{Definition}
\newtheorem{defn}[equation]{Definition}

\newtheorem{rem}[equation]{Remark}

\newfont{\bigbf}{cmbx10 scaled\magstep1}
\normalbaselines
\numberwithin{equation}{section}

\def\no{\noindent}

\def\C{{\mathbb C}}
\def\R{{\mathbb R}}
\def\H{{\mathbb H}}

\def\N{{\mathbb N}}

\def\Z{{\mathbb Z}}

\def\al{\alpha}

\def\ga{\gamma}
\def\Ga{\Gamma}
\def\de{\delta}
\def\De{\Delta}
\def\eps{\epsilon}
\def\la{\lambda}
\def\La{\Lambda}
\def\si{\sigma}

\def\Om{\Omega}

\def\3{\ss}

\def\acts{\curvearrowright}

\def\D{\partial}

\def\diamo{\diamondsuit}
\def\embed{\hookrightarrow}

\def\Flag{\operatorname{Flag}}

\def\Flagn{\Flag(\nu_{mod})}
\def\Flagt{\Flag(\tau_{mod})}

\def\geo{\partial_{\infty}}

\def\Hom{\operatorname{Hom}}

\def\interior{\operatorname{int}}
\def\Isom{\mathop{\hbox{Isom}}}

\def\lra{\longrightarrow}

\def\midp{\operatorname{mid}}

\def\pihalf{\frac{\pi}{2}}

\def\2pithird{\frac{2\pi}{3}}

\def\st{\operatorname{st}}
\def\ost{\mathop{\hbox{ost}}}

\def\tangle{\angle_{Tits}}
\def\tauto{\buildrel f\over\lra}
\def\tits{\partial_{Tits}}
\def\Th{\mathop{\hbox{Th}}\nolimits}

\def\diamot{\diamondsuit_{\tau_{mod}}}
\def\diamoTh{\diamondsuit_{\Theta}}
\def\diamoThp{\diamondsuit_{\Theta'}}

\def\Lat{\La_{\tau_{mod}}}
\def\inte{\operatorname{int}}

\def\interior{\operatorname{int}}

\def\DF{\partial_{F\ddot u}}

\def\ThF{\operatorname{Th}_{F\ddot u}}

\def\F{{\mathrm F}}
\def\8{\infty}
\def\<{\langle}
\def\>{\rangle}

\def\BI{\begin{itemize}}
\def\EI{\end{itemize}}

\def\numod{\nu_{mod}}

\def\taumod{\tau_{mod}}
\def\simod{\sigma_{mod}}

\title{Morse actions of discrete groups on symmetric spaces: Local-to-global principle}
\author{Michael Kapovich, Bernhard Leeb, Joan Porti}
\date{\today}

\begin{document}

\maketitle


\begin{abstract}
\no
Our main result is a local-to-global principle 
for Morse quasigeodesics, maps and actions. 
As an application of our techniques 
we show algorithmic recognizability 
of Morse actions and construct Morse ``Schottky subgroups'' 
of higher rank semisimple Lie groups via arguments 
not based on Tits' ping-pong.
Our argument is purely geometric 
and proceeds by constructing    
equivariant Morse quasiisometric embeddings of trees 
into higher rank symmetric spaces. 
\end{abstract}

\tableofcontents

\section{Introduction}

This is a sequel to our paper \cite{anolec} and mostly consists of the material of section 7 of our earlier paper \cite{morse} (the only additional material appears in Theorem \ref{thm:Omstability} 
and the appendix to the paper). We recall that 
quasigeodesics in Gromov hyperbolic spaces can be recognized locally 
by looking at sufficiently large finite pieces, see \cite{CDP}. In our earlier papers \cite{coco15, anolec, mlem, bordif, manicures}, for higher rank symmetric spaces $X$ (of noncompact type)  we introduced an analogue of hyperbolic quasigeodesics, which we call {\em Morse quasigeodesics}. Morse quasigeodesics are defined relatively to a certain face $\taumod$ 
of the model spherical face $\simod$ of $X$. In addition to the quasiisometry constants $L, A$, $\taumod$-Morse   quasigeodesics come equipped with two other parameters, a positive number $D$ and a {\em Weyl-convex} subset $\Theta$ of the {\em open star} of $\taumod$ in the modal spherical chamber $\simod$.   In \cite{morse, anolec, mlem} we also defined $\taumod$-Morse maps $Y\to X$ 
from Gromov-hyperbolic spaces to symmetric spaces. These maps are defined by the property that they send geodesics to {\em uniformly $\taumod$-Morse quasigeodesics}, i.e. $\taumod$-Morse quasigeodesics with a fixed set of parameters, $(\Theta, D, L, A)$.

The main result of this paper is a 
{\em local} characterization of Morse quasigeodesics in $X$: 


\begin{thm}[Local-to-global principle for Morse quasigeodesics]
\label{thm:locglobmqg0}
For $L,A,\Theta,\Theta',D$ there exist $S,L',A',D'$ such that 
every $S$-local $(\Theta, D, L, A)$-local Morse quasigeo\-de\-sic in $X$ 
is a $(\Theta', D', L', A')$--Morse quasigeodesic. 
\end{thm}

Here $S$-locality of a certain property of a map  means that this property is satisfied for restrictions of this map to subintervals of 
length $S$. We refer  to Definition \ref{thm:locglobmqg} and Theorem \ref{thm:locglobmqg} for the details. Based on this principle, we 
prove in Section \ref{sec:locglobmm} a local-to-global principle for Morse maps from hyperbolic metric spaces to symmetric spaces.

\medskip 
We  prove several consequences of these local-to-global principles:

1. The {\em structural stability} of Morse  subgroups of $G$, 
generalizing Sullivan's Structural Stability Theorem in rank one \cite{Sullivan} (see also \cite{KKL} for a detailed proof);
see Theorems~\ref{thm:morsestab} and~\ref{thm:strcstab}. 
While structural stability for Anosov subgroups was known earlier (Labourie and Guichard--Wienhard), 
our method is more general and applies to a wider class of discrete subgroups, see \cite{relmorse-2}.

\begin{thm}
[Openness of the space of Morse actions] 
For a word hyperbolic group $\Ga$, the subset of $\tau_{mod}$-Morse actions is open in $\Hom(\Ga, G)$. 
\end{thm}

\begin{thm}
[Structural stability] \label{thm:stability}
Let $\Ga$ be word hyperbolic.
Then for $\tau_{mod}$-Morse actions $\rho:\Ga\acts X$, 
the boundary embedding $\al_{\rho}: \geo \Ga\to \Flagt$ depends continuously on the action $\rho$. 
\end{thm}

In particular, actions sufficiently close to a faithful Morse action 
are again discrete and faithful. We supplement this structural stability theorem with a stability theorem on {\em domains of proper discontinuity}, Theorem \ref{thm:Omstability}.

2. The locality of the Morse property implies  that Morse subgroups are algorithmically recognizable; Section \ref{sec:algrec}:

\begin{thm}
[Semidecidability of Morse property of group actions] 
Let $\Ga$ be word hyperbolic. 
Then there exists an algorithm whose inputs are homomorphisms $\rho: \Ga\to G$ (defined on generators of $\Ga$) 
and which terminates if and only if $\rho$ defines a $\tau_{mod}$-Morse action $\Ga\acts X$. 
\end{thm}

If the action is not Morse, the algorithm runs forever. Note that in view of \cite{Kapovich2015}, there are no algorithms (in the sense of BSS computability) which would recognize if a representation $\Ga\to \Isom(\H^3)$ is {\em not geometrically finite}.

\medskip 
3. We illustrate our techniques by constructing Morse-Schottky actions of free groups on higher rank symmetric spaces; Section \ref{sec:schottky actions}.
Unlike all previously known constructions, our proof does not rely on ping-pong arguments, but
is purely geometric and proceeds by constructing  equivariant quasi-isometric embeddings of trees. 
The key step is the observation that a certain local {\em straightness} property for 
sufficiently spaced sequences of points in the symmetric space implies the global Morse property. 
This observation is also at the heart of the proof 
of the local-to-global principle for Morse actions.

\medskip
Since \cite{morse} was originally posted in 2014, several improvements on the material of section 7 of \cite{morse} and, hence, of the present paper were made:

\medskip 
(a) Different forms of Combination Theorems for Anosov subgroups were proven in \cite{DKL, DK1, DK2}  
written in collaboration with Subhadip Dey by the 1st and the 2nd author and, subsequently, by the 1st author. The first one was a generalization of the technique in section \ref{sec:schottky actions} the present paper, but the other two generalizations are based on a form of the ping-pong argument. 

\medskip 
(b) Explicit estimates in the local-to-global principle for Morse quasigeodesics and, hence, Morse embeddings, were obtained by Max Riestenberg in  \cite{Riestenberg}. Riestenberg's estimates are based on replacing certain limiting arguments used in the present paper with differential-geometric and Lie-theoretic arguments.

\newpage
{\bf Organization of the paper.} 

The notions of Morse quasigeodesics and actions 
are discussed in detail in section \ref{sec:morse}. 
In that section, among other things, 
we establish local-to-global principles for Morse quasigeodesics.


In section \ref{sec:group-applications} we apply local-to-global principles to discrete subgroups of Lie groups: We show that Morse actions are structurally  stable and algorithmically recognizable. We also construct 
Morse-Schottky actions of free groups on symmetric spaces.    In section 5 (the appendix to the paper) we prove further properties 
of Morse quasigeodesics that we found to be useful in our work.

\medskip 
{\bf Acknowledgements.} 
The first author was supported by NSF grants  DMS-12-05312 and DMS-16-04241, by 
KIAS (the Korea Institute for Advanced Study) through the KIAS scholar program, 
and by a Simons Foundation Fellowship, grant number 391602. 
The last author was supported by grants  Mineco MTM2012-34834 and AGAUR  SGR2009-1207. The three authors are also grateful to 
the GEAR grant which partially supported the IHP trimester in Winter of 2012 (DMS 1107452,
1107263, 1107367 ``RNMS: Geometric structures and representation varieties'' 
(the GEAR Network), and to the Max Planck Institute for Mathematics in Bonn, where some of this work was done. 

\section{Preliminaries}

\subsection{Basic notions of geometry of symmetric spaces}  

Throughout the paper we will be using definitions, notations and results of our earlier work.

We refer the reader to our earlier papers, e.g. \cite{coco15, anolec, mlem, bordif, manicures} for the various notions related to symmetric spaces, 
such as {\em polyhedral Finsler metrics} on symmetric spaces (\cite{bordif}),  the {\em opposition involution} $\iota$ of $\simod$, {\em model faces} $\taumod$ of $\simod$ and the associated $\taumod$-flag manifolds $\Flagt$  
(sections 2.2.2 and 2.2.3 of \cite{anolec}), {\em type map} $\theta: \geo X\to \simod$, 
 open Schubert cells $C(\tau)\subset \Flagt$ (section 2.4 of \cite{anolec}), 
$\Delta$-valued distances $d_\Delta$  on $X$ (section 2.6 of \cite{anolec}), 
 $\Theta$-regular geodesic segments (see \S 2.5.3 of \cite{anolec}),  
{\em parallel sets}, {\em stars}, {\em open stars} 
and $\Theta$-stars, 
$\st(\tau)$, $\ost(\tau)$, and $\st_\Theta(\tau)$, {\em Weyl sectors} $V(x, \tau)$ (section 2.4 of \cite{anolec}),  
 {\em Weyl cones} $V(x, \st(\tau))$ and $\Theta$-cones $V(x, \st_\Theta(\tau))$, 
{\em diamonds} $\diamot(x,y)$ and $\Theta$-diamonds $\diamoTh(x,y)$  (section 2.5 of \cite{anolec}), {\em $\taumod$-regular sequences and groups} (section 4.2 of \cite{anolec}), $\taumod$-convergence subgroups, flag-convergence,  the Finsler interpretation of 
flag-convergence (see \cite[\S 4.5 and 5.2]{bordif} and \cite{anolec}), 
$\taumod$-limit sets $\Lat(\Ga)\subset \Flagt$ (section 4.5 of \cite{anolec}), visual limit set (page 4 of \cite{anolec}), uniformly $\taumod$-regular sequences and subgroups (section 4.6 of \cite{anolec}), Morse 
subgroups (section 5.4 of \cite{anolec}) and, more generally, Morse quasigeodesics and Morse maps (Definitions 5.31, 5.33 of \cite{mlem}), antipodal limit sets (Definition. 5.1 of \cite{anolec}) and antipodal maps to flag-manifolds (Definition 6.11 of \cite{mlem}).

\medskip 
In the paper we will be frequently using convexity of $\Theta$-cones in $X$:

\begin{prop}[Proposition 2.10  in \cite{anolec}]
\label{prop:thconeconv}
For every Weyl-convex subset $\Theta\subset \st(\taumod)$, for every $x\in X$ and $\tau\in \Flagt$, the cone $V(x, \st_\Theta(\tau))\subset X$ 
is convex. 
\end{prop}

\subsection{Standing notation and conventions} 

\begin{itemize}

\item We will use the notation $X$ for a symmetric space of {\em noncompact type}, $G$ for 
 a semisimple Lie group 
acting isometrically and transitively on $X$, 
and $K$ is a maximal compact subgroup of $G$, so that $X$ is diffeomorphic to $G/K$.   
We will assume that $G$ is commensurable with the isometry group $\Isom(X)$ 
in the sense that we allow finite kernel and cokernel 
for the natural map $G\to\Isom(X)$. In particular, 
the image of $G$ in $\Isom(X)$ contains the identity component $\Isom(X)_o$. 

\item We let $\tau_{mod}\subseteq\si_{mod}$ be a fixed  $\iota$-invariant face type.

\item We will use the notation $x_n\tauto \tau\in \Flagt$ for the {\em flag-convergence} of a $\taumod$-regular sequence $x_n\in X$ to a simplex $\tau\in \Flagt$. 

\item We will be using the notation $\Theta, \Theta'$ for an $\iota$-invariant,  
compact, Weyl-convex (see Definition 2.7 in \cite{anolec}) 
subset of the open star $\ost(\taumod)\subset \simod$. 

\item We will always assume that $\Theta< \Theta'$, meaning that $\Theta\subset \interior(\Theta')$.

\item  Constants $L, A, D, \eps, \de, l, a, s, S$ 
are meant to be always strictly positive and $L\ge 1$. 
\end{itemize}

\subsection{$\zeta$-angles}\label{sec:zangle}

\def\zangle{\angle^{\zeta}}

We fix as auxiliary datum 
a $\iota$-invariant type $\zeta=\zeta_{mod}\in\interior(\tau_{mod})$. (We will omit the subscript in $\zeta_{mod}$ 
in order to avoid cumbersome notation for $\zeta$-angles.)  For a simplex $\tau\subset\geo X$ of type $\tau_{mod}$, i.e. $\tau\in \Flagt$,  
we define $\zeta(\tau)\in\tau$ as the ideal point of type $\zeta_{mod}$.  Given two such simplices 
$\tau_\pm\in \Flagt$ and a point $x\in X$,  define 
the {\em $\zeta$-angles} 
\begin{equation}\label{eq:zeta-angle}
\zangle_x (\tau_-, \tau_+)= \zangle_x (\tau_-, \xi_+):=  
\angle_x(\xi_-, \xi_+), 
\end{equation}
where $\xi_\pm =\zeta(\tau_\pm)$.

 Similarly, define the {\em $\zeta$-Tits angle} 
\begin{equation}\label{eq:zeta-tangle}
\tangle^\zeta(\tau_-, \tau_+)= \tangle^\zeta(\tau_-, \xi_+):= 
\angle_x(\xi_-, \xi_+),   
\end{equation}
where $x$ belongs to a flat $F\subset X$ such that $\tau_-, \tau_+\subset \tits F$. Then simplices $\tau_\pm$ 
(of the same type) are antipodal iff
$$
\tangle^\zeta(\tau_-, \tau_+)=\pi
$$
for some, equivalently, every, choice of $\zeta$ as above.

\begin{rem}
\label{rem:antip}
We observe that the ideal points $\zeta_{\pm}$ are opposite,
$\tangle(\zeta_-,\zeta_+)=\pi$, 
if and only if they can be seen under angle $\simeq\pi$ (i.e., close to $\pi$) 
from some point in $X$. More precisely, 
there exists $\eps(\zeta_{mod})$ such that:
 
{\em If $\angle_x(\zeta_-,\zeta_+)>\pi-\eps(\zeta_{mod})$ 
for some point $x$ then $\zeta_{\pm}$ are opposite. }

\noindent This follows from the angle comparison 
$\angle_x(\zeta_-,\zeta_+)\leq\tangle(\zeta_-,\zeta_+)$
and the fact that the Tits distance between ideal points 
of the fixed type $\zeta_{mod}$ 
takes only finitely many values. 
\end{rem}

For a $\tau_{mod}$-regular unit tangent vector $v\in TX$ 
we denote by $\tau(v)\subset\geo X$ the unique simplex of type $\tau_{mod}$ 
such that ray $\rho_v$ with  the initial direction $v$ 
represents an ideal point in $\ost(\tau(v))$. 
We put $\zeta(v)=\zeta(\tau(v))$. 
Note that $\zeta(v)$ depends continuously on $v$.

For a $\taumod$-regular segment $xy$ in $X$ we let $\tau(xy)=\tau(v)$,
where $v$ is the unit vector tangent to $xy$.

Then, for a $\taumod$-regular segments $xy, xz$ and $\tau\in \Flagt$, 
we define the $\zeta$-angles
$$
\zangle_x(y, \tau) = \zangle_x(\tau(xy), \tau), \quad \zangle_x(y, z) = \zangle_x(\tau(xy), \tau(xz)) 
$$
 
Thus, the $\zeta$-angle depends not on $y, z$ but rather  on the simplices $\tau(xy), \tau(xz)$. 
These $\zeta$-angles  will play the role of angles the between diamonds 
$\diamot(x, y)$ and $\diamot(x, z)$, meeting  at $x$. Note that if $X$ has rank 1, then the 
$\zeta$-angles are just the ordinary Riemannian angles.

\subsection{Distances to parallel sets versus angles}
\label{sec:angledist}

In this section we collect some  geometric facts regarding 
parallel sets in symmetric spaces, primarily dealing with estimation of distances from points in $X$ to parallel sets. 

\begin{rem}
The constants and functions in this section are not explicit and their existence is proven by compactness arguments. 
For explicit computations here and in Theorem \ref{thm:locstrimplcoastrseq}, we refer the reader to the PhD thesis of ... 
\end{rem}

 We first prove a lemma (Lemma \ref{lem:strongly asymptotic geodesic to a star}) 
 which strengthens Corollary 2.46 of \cite{anolec}.


\begin{lemma}\label{lem:strongly asymptotic geodesic to a star} 
Suppose that $\tau_\pm$ are antipodal simplices in $\tits X$. Then 
every geodesic ray $\ga$ asymptotic to a point $\xi\in \ost(\tau_+)$, 
is strongly asymptotic to a geodesic ray in $P(\tau_-, \tau_+)$. 
\end{lemma} 
\proof If $\xi$ belongs to the interior of the simplex $\tau_+$, then the assertion follows from 
Corollary 2.46 of \cite{anolec}: 

\medskip 
{\em Weyl sectors $V(x_1,\tau)$ and $V(x_2,\tau)$ are strongly asymptotic 
if and only if $x_1$ and $x_2$ lie in the same horocycle at $\tau$. }

\medskip 
We now consider the general case. Suppose, that $\xi$ belongs to an open simplex 
$\interior(\tau')$, such that $\tau$ is a face of $\tau'$. Then there exists an apartment $a\subset \tits X$ containing both $\xi$ (and, hence, $\tau'$ as well as $\tau$) and the simplex $\tau_-$. Let $F\subset X$ be the maximal flat with $\geo F=a$. Then $F$ contains a geodesic asymptotic to points in 
$\tau_-$ and $\tau_+$. Therefore, $F$ is contained in $P(\tau_-, \tau_+)$. On the other hand, by the same
Corollary 2.46 of \cite{anolec}, 
applied to the simplex $\tau'$, we conclude that  $\ga$ is strongly asymptotic to a geodesic ray in $F$. \qed 

\medskip 
The following lemma provides a quantitative strengthening of the conclusion of Lemma \ref{lem:strongly asymptotic geodesic to a star}: 

\begin{lem}
\label{lem:decay} 
Let $\Theta$ be a compact subset of $\ost(\tau_+)$. Then 
those rays $x\xi$ with $\theta(\xi)\in\Theta$ 
are uniformly strongly asymptotic to $P(\tau_-,\tau_+)$, 
i.e.\ $d(\cdot,P(\tau_-,\tau_+))$ decays to zero along them 
uniformly in terms of $d(x, P(\tau_-, \tau_+))$ and $\Theta$. 
\end{lem}
\proof Suppose that the assertion of lemma is false, i.e., there exists $\eps>0$, a sequence $T_i\in \R_+$ diverging to infinity, and 
a sequence of  rays  $\rho_i=x_i \xi_i$ with  $\xi_i\in \Theta$ and $d(x_i, P(\tau_-,\tau_+))\le d$, so that
\begin{equation}\label{eq:outside}
d(y, P(\tau_-,\tau_+))\ge \eps, \forall y\in \rho([0, T_i]).
\end{equation}
Using the action of the stabilizer of 
$P(\tau_-,\tau_+)$, we can assume that the points $x_i$ belong to a certain compact subset of $X$. Therefore, the sequence of rays 
$x_i \xi_i$ subconverges to a ray $x\xi$ with $d(x, P(\tau_-,\tau_+))\le d$ and $\xi\in \Theta$. The inequality \eqref{eq:outside} then implies that 
the entire limit ray $x\xi$ is contained outside of the open $\eps$-neighborhood of the parallel set $P(\tau_-,\tau_+)$. However, in view of Lemma 
\ref{lem:strongly asymptotic geodesic to a star}, the ray $x\xi$ is strongly asymptotic to a geodesic in $P(\tau_-,\tau_+)$. Contradiction. \qed

\medskip

We next relate distances from points $x\in X$ to parallel sets and the $\zeta$-angles at $x$.  Suppose that the simplices $\tau_{\pm}$, equivalently, 
the ideal points $\zeta_{\pm}=\zeta(\tau_\pm)$ (see section \ref{sec:zangle}), are opposite. 
Then 
$$
\zangle_x(\tau_-, \tau_+)= \angle_x(\zeta_-,\zeta_+)=\pi$$  
if and only if $x$ lies in the parallel set $P(\tau_-,\tau_+)$.  
Furthermore, 
$\zangle_x(\tau_-,\tau_+)\simeq\pi$
if and only if $x$ is close to $P(\tau_-,\tau_+)$, 
and both quantities control each other near the parallel set. 
More precisely: 

\begin{lem}
\label{lem:distangcontr}
(i)
If $d(x,P(\tau_-,\tau_+))\leq d$, 
then $\zangle_x(\tau_-,\tau_+)\geq\pi-\eps(d)$ 
with $\eps(d)\to0$ as $d\to0$. 

(ii)
For sufficiently small $\eps$, $\eps\leq\eps'(\zeta_{mod})$, we have:
The inequality $\zangle_x(\tau_-,\tau_+)\geq\pi-\eps$ implies that  
$d(x,P(\tau_-,\tau_+))\leq d(\eps)$ for some function $d(\eps)$ which converges to $0$ as $\eps\to0$.
\end{lem}
\proof
The intersection of parabolic subgroups 
$P_{\tau_-}\cap P_{\tau_+}$
preserves the parallel set $P(\tau_-,\tau_+)$ and acts transitively on it.
Compactness and the continuity of $\angle_{\cdot}(\zeta_-,\zeta_+)$
therefore imply that $\pi-\angle_{\cdot}(\zeta_-,\zeta_+)$ 
attains on the boundary of the tubular $r$-neighborhood of $P(\tau_-,\tau_+)$ 
a strictly positive maximum and minimum, 
which we denote by $\phi_1(r)$ and $\phi_2(r)$. 
Furthermore, $\phi_i(r)\to0$ as $r\to0$. 
We have the estimate:
\begin{equation*}
\pi-\phi_1(d(x,P(\tau_-,\tau_+)))
\leq\angle_x(\zeta_-,\zeta_+)\leq
\pi-\phi_2(d(x,P(\tau_-,\tau_+)))
\end{equation*}
The functions $\phi_i(r)$ are (weakly) monotonically increasing. 
This follows from the fact that, 
along rays asymptotic to $\zeta_-$ or $\zeta_+$, 
the angle $\angle_{\cdot}(\zeta_-,\zeta_+)$ is monotonically increasing 
and the distance $d(\cdot,P(\tau_-,\tau_+))$ is monotonically decreasing. 
The estimate implies the assertions. 
\qed

\medskip
The control of 
$d(\cdot,P(\tau_-,\tau_+))$  and $\angle_{\cdot}(\zeta_-,\zeta_+)$ 
``spreads'' along the Weyl cone $V(x,\st(\tau_+))$, 
since the latter is asymptotic to the parallel set $P(\tau_-,\tau_+)$. 
Moreover, the control improves,  
if one enters the cone far into a $\tau_{mod}$-regular direction. 
More precisely: 
\begin{lem}
\label{lem:distangdec}
Let $y\in V(x,\st_{\Theta}(\tau_+))$
be a point with $d(x,y)\geq l$. 

(i) 
If $d(x,P(\tau_-,\tau_+))\leq d$,  
then 
\begin{equation*}
d(y,P(\tau_-,\tau_+))\leq D'(d,\Theta,l)\leq d
\end{equation*}
with $D'(d,\Theta,l)\to0$ as $l\to+\infty$. 

(ii)
For sufficiently small $\eps$, $\eps\leq\eps'(\zeta_{mod})$, we have:
If $\angle_x(\zeta_-,\zeta_+)\geq\pi-\eps$, 
then 
\begin{equation*}
\angle_y(\zeta_-,\zeta_+)\geq\pi-\eps'(\eps,\Theta,l)\geq\pi-\eps(d(\eps))
\end{equation*}
with $\eps'(\eps,\Theta,l)\to0$ as $l\to+\infty$. 
\end{lem}
\proof
The distance from $P(\tau_-,\tau_+)$
takes its maximum at the tip $x$ of the cone $V(x,\st(\tau_+))$, 
because it is monotonically decreasing along the rays $x\xi$ 
for $\xi\in\st(\tau_+)$. 
This yields the right-hand bounds $d$ and, 
applying Lemma~\ref{lem:distangcontr} twice, $\eps(d(\eps))$. 

Those rays $x\xi$ 
with uniformly $\tau_{mod}$-regular type $\theta(\xi)\in\Theta$ 
are uniformly strongly asymptotic to $P(\tau_-,\tau_+)$, 
i.e.\ $d(\cdot,P(\tau_-,\tau_+))$ decays to zero along them 
uniformly in terms of $d$ and $\Theta$, see Lemma \ref{lem:decay}.  
This yields the decay $D'(d,\Theta,l)\to0$ as $l\to+\infty$. 
The decay of $\eps'$ follows by applying Lemma~\ref{lem:distangcontr} again. 
\qed

\section{Morse maps}
\label{sec:morse}

In this section we  investigate the Morse property of sequences and maps. 
The main aim of this section is to establish 
a local criterion for being Morse.
To do so we introduce a local notion of  \emph{straightness} for sequences of points in $X$. 
Morse sequences are in general not straight, but they become straight after suitable modification, namely by sufficiently
coarsifying them
and then passing to the sequence of successive midpoints.
Conversely, the key result is that sufficiently spaced  straight sequences are Morse. We conclude that there is a local-to-global characterization of the Morse property.

\subsection{Morse quasigeodesics}

\begin{dfn}[Morse quasigeodesic]
\label{dfn:mqg}
A $(\Theta,D,L,A)$-{\em Morse quasigeodesic} in $X$ 
is an $(L,A)$-quasigeodesic $p: I \to X$ (defined on an interval $I\subset \R$) 
such that for all $t_1,t_2\in I$  
the subpath $p|_{[t_1,t_2]}$ is $D$-close to a $\Theta$-diamond 
$\diamoTh(x_1,x_2)$ 
with $d(x_i,p(t_i))\leq D$.
\end{dfn}

We will refer to a quadruple $(\Theta,D,L,A)$ as a {\em Morse datum}  and abbreviate $M=(\Theta,D,L,A)$. Set  $M+D'=(\Theta,D+D',L,A +2D')$. We say that 
$M$ contains $\Theta$  if $M$ has the form 
$(\Theta, D, L, A)$ for some $D\ge 0, L\ge 1, A\ge 0$.


\medskip 
The following lemma is immediate from the definiton of a $M$-Morse quasigeodesic.

\begin{lemma}
[Perturbation lemma]\label{lem:perturbation}
If $p, p'$ are paths in $X$ such that $p$ is $M$-Morse and $d(p, p')\le D'$ then $p'$ is $M+D'$-Morse. 
\end{lemma}

\medskip
A Morse quasigeodesic $p$ is called a {\em Morse ray} if its domain is a half-line. If $I=\R$ then a Morse quasigeodesic is called a {\em Morse quasiline}. 

\medskip 
Morse quasirays do in general not converge at infinity (in the visual compactification of $X$),  
but they $\tau_{mod}$-converge at infinity. 
This is a consequence of:

\begin{lem}[Conicality]
\label{lem:morsecon}
Every Morse quasiray $p: [0,\infty)\to X$ 
is uniformly Hausdorff close to a subset of a cone $V(p(0),\st_{\Theta}(\tau))$
for a unique simplex $\tau$ of type $\tau_{mod}$. 
\end{lem}
\proof
The subpaths $p|_{[0,t_0]}$ are uniformly Hausdorff close 
to $\Theta$-diamonds. 
These subconverge to a cone $V(x,\st_{\Theta}(\tau))$ 
$x$ uniformly close to $p(0)$ and $\tau$ a simplex of type $\tau_{mod}$. 
This establishes the existence. 
Since $p(n)\tauto\tau$, 
the uniqueness of $\tau$ follows from the uniqueness of $\tau_{mod}$-limits, 
 see \cite[Lemma 4.23]{anolec}. 
\qed

\begin{dfn}[End of Morse quasiray]
We call the unique simplex given by the previous lemma 
the {\em end} of the Morse quasiray $p:[0,\infty)\to X$ 
and denote it by 
\begin{equation*}
p(+\infty)\in\Flag(\tau_{mod}).
\end{equation*}
\end{dfn}
Hausdorff close Morse quasirays have the same end by Lemma \ref{lem:morsecon}.  In section \ref{sec:continuity} we will 
prove uniform continuity of ends of Morse quasirays with respect to the topology of {\em coarse convergence} of quasirays.



\subsection{Morse maps}

We now turn to {\em Morse maps} with more general domains (than just intervals). 

\begin{defn}
Let $Y$ be a Gromov-hyperbolic geodesic metric space. A map $f: Y\to X$ is called 
$M$-Morse if it sends geodesics in $Y$ to $M$-Morse quasigeodesics. 
\end{defn}

Thus, every Morse map is a quasiisometric embedding. While this definition makes sense for general metric spaces,
 in \cite{mlem} we proved that the domain of a Morse map is necessarily hyperbolic. 

\medskip 
More generally, one can define Morse maps on {\em quasigeodesic metric spaces}:

\begin{dfn}[Quasigeodesic metric space]
A metric space $Z$ is called {\em $(l, a)$-quasigeodesic} 
if all pairs of points in $Y$ can be connected by $(l, a)$-quasigeodesics.  A space  is called {\em quasigeodesic} 
if it is $(l, a)$-quasigeodesic for some pair of parameters $l, a$.
\end{dfn}

Every quasigeodesic space  is quasiisometric to a geodesic metric space. Namely, if $Z$  is 
$(\la, \al)$-quasigeodesic space then it is quasiisometric to its $(\la+\al)$-Rips complex.  
The quasigeodesic spaces considered in this paper 
are discrete groups equipped with word metrics.

\begin{dfn}[Morse embedding]
\label{dfn:morseemb}
Let $(\Theta,D,L,A)$ be a Morse datum. \newline 
An {\em $(\Theta,D,L,A,l,a)$-Morse embedding} (or a map) 
from an $(l,a)$-quasigeodesic space $Z$ into $X$ is a map 
$f: Z\to X$ which sends  $(l,a)$-quasigeodesics in $Z$ 
to $(\Theta,D,L,A)$-Morse quasigeodesics in $X$.  
\end{dfn}

Of course, every $(l,a)$-quasigeodesic metric space is also $(l',a')$-quasigeodesic space for any $l'\ge l, a'\ge a$. 
The next lemma shows that this choice of quasigeodesic constants is essentially irrelevant.

\begin{lem}
\label{lem:morsefromhyp}
Let $f:Z\to X$ be a map from a  Gromov-hyperbolic $(l,a)$-quasigeodesic space $Z$. 
If $f$ is $M=(\Theta,D,L,A,l,a)$-Morse then for any $(l',a')$, it sends  $(l',a')$-quasigeodesics in $Z$ to 
$M'=(\Theta,D',L',A')$-Morse quasigeodesics in $X$. Here the datum $M'$ depends only on $M, l', a'$ and 
the hyperbolicity constant $\delta$ of $Z$. 
\end{lem}
\proof
This is a consequence of the definition of Morse quasigeodesics,  
and the Morse Lemma applied to $Z$.  
\qed
 
Notice that  the parameter $\Theta$ in the Morse datum $M'$ is the same as in $M$. Hence, we 
arrive to

\begin{definition}
A map $f: Z\to X$ of a quasigeodesic hyperbolic space $Z$ is called $\Theta$-{\em Morse} if it sends  uniform quasigeodesics  
in $Z$ to $\Theta$-Morse uniform quasigeodesics in $X$.  
\end{definition} 

This notion depends only on the quasi-isometry class of $Z$, 
i.e.\  
the precomposition of a $\Theta$-Morse embedding with a quasi-isometry 
is again $\Theta$-Morse. 
For this to be true  we have to require control on the images of quasigeodesics 
of arbitrarily bad (but uniform) quality.

\medskip

Let $\Ga$ be a hyperbolic group with fixed a finite generating set $S$, and let $Y$ be the  
Cayley graph of $\Ga$ with respect to $S$. For $x\in X$, an isometric action $\Ga\acts X$ determines the {\em orbit map} $o_x: \Ga\to \Ga x\subset X$. Every such map extends to the Cayley graph $Y$ of 
$\Ga$, sending edges to geodesics in $X$.

\begin{defn}
An isometric action $\Ga\acts X$ or a representation $\rho: \Ga\to G$,  
 is called $M$-Morse (with respect to a base-point $x\in X$) if the (extended) orbit map 
$o_x: Y\to X$ is $M$-Morse. Similarly, a subgroup $\Ga< G$ is  Morse if the inclusion homomorphism 
 $\Ga\embed G$ is Morse.  
\end{defn}

The Morse property of an action and the parameter $\Theta$, of course, does not depend on the choice of a generating set of 
$\Ga$ and a  base-point $x$, but the triple  $(D,L,A)$ does.  Thus, it makes sense to talk about a $\Theta$-Morse and $\taumod$-Morse actions of hyperbolic groups, where $\Theta\subset \ost(\taumod)$. 
In \cite{anolec, mlem, bordif} we gave many alternative definitions of Morse actions, including the equivalence of this definition to the notion of Anosov subgroups.

\subsection{Continuity at infinity}\label{sec:continuity}

Let $X, Y$ be proper metric spaces. We fix a base point $y\in Y$.

\begin{defn}\label{def:coconv}
A sequence of maps $f_n: Y\to X$ is said to {\em coarsely  converge} 
to a map $f: Y\to X$ if there exists $C<\infty$ such that 
for every $R$  there exists $N=N(C,R)$ for which
$$
d(f_n|_B, f|_B)\le C, 
$$ 
where $B=B(y,R)$. 
\end{defn}

Note the difference of this definition with the notion of uniform convergence on compacts: Since we are working in the coarse setting, 
requiring the distance between maps to be less than $\eps$ close to zero is pointless. 

In view of the Arzela--Ascoli theorem, the space of $(L,A)$-coarse Lipschitz maps $Y\to X$ sending $y$ to a fixed bounded subset of $X$, 
is coarsely sequentially compact: Every sequence contains a coarsely converging subsequence.

\medskip 
In the next lemma we assume that $Y$ is a geodesic $\delta$-hyperbolic space and $X$ is a symmetric space of noncompact type. The lemma itself is an immediate consequence 
of the perturbation lemma, Lemma \ref{lem:perturbation}.

\begin{lemma}
Suppose that $p_n: \R_+\to X$ is a sequence of $M$-Morse rays which coarsely converges to a map $p: \R_+\to X$. 
Then $p$ is $M'$-Morse, where $M'=M+C$ and the constant $C$ is the one appearing in the definition of coarse convergence. 
\end{lemma}

In particular, a coarse limit of a sequence of (uniformly) Morse quasigeodesics is again Morse. 

\medskip

For the next lemma, we equip the flag manifold $\F=\Flagt$ with some background metric $d_\F$. 

\begin{lemma}\label{lem:cont-dependence}
Suppose that $p_n :\R_+\to X$ is a sequence of $M$-Morse rays coarsely converging to a $M$-Morse ray  
$p: \R_+\to X$. Then the sequence $\tau_n:= p_n(\infty)$ of ends of the quasirays $p_n$ 
converges to $\tau=p(\infty)$. Moreover, the latter convergence is uniform in the following sense. 
For every $\eps>0$ there exists $n_0$ depending only on $M$ and $C$ and $N(R,C)$  
(appearing in Definition \ref{def:coconv}) such that  for all $n\ge n_0$, $d_{\F}(\tau_n, \tau)\le \eps$. 
\end{lemma}
\proof Suppose that the claim is false. Then in view of coarse compactness of the space of $M$-Morse maps sending $y$ to 
a fixed compact subset of $X$, there exists a sequence $(p_n)$  as in the lemma, coarsely converging to $p$, such that the sequence 
$p_n(\infty)=\tau_n$ converges to $\tau'\ne p(\infty)=\tau$. By  the coarse convergence $p_n\to p$, 
there exists $C< \infty$ and a sequence $t_n\to \infty$  such that $d(p_n(t_n), p(t_n))\le C$. By the definition of Morse quasigeodesics, 
there exists a sequence of cones $V(x_n, \st(\tau_n))$ (with $x_n$ in a bounded subset $B\subset X$) 
such that the image of $p_n$ is contained in the $D$-neighborhood of $V(x_n, \st(\tau_n))$. Thus, the sequence $(p_n(t_n))$ flag-converges to 
$\tau'$, while $(p(t_n))$ flag-converges to $\tau$. According to \cite[Lemma 4.23]{anolec}, altering a sequence by a uniformly bounded amount, does not change the flag-limit. Therefore, the sequence $(p(t_n))$ also flag-converges to $\tau'$. Hence, $\tau=\tau'$. A contradiction. \qed

\subsection{A Morse Lemma for straight sequences}

In order to motivate the results of this section we recall the following {\em sufficient condition} for a piecewise-geodesic path in a Hadamard manifold $Y$ of curvature $\le -1$ to be quasigeodesic (see e.g. \cite{KaLi}):

\begin{prop}
Suppose that $c$ is a piecewise-geodesic path in $Y$ whose angles at the vertices are $\ge \al>0$ and whose edges are longer than $L$, where $\al$ and $L$ satisfy
\begin{equation}\label{eq:rank1-qg}
\cosh(L /2)\sin(\alpha/2) \ge \nu>1.  
\end{equation}
Then $c$ is an $(L(\nu), A(\nu))$-quasigeodesic. \end{prop}

By considering $c$ with vertices on a horocycle in the hyperbolic plane, one see that the inequality in this proposition is sharp.

\begin{cor}\label{cor:rank1} 
If $L$ is sufficiently large and $\alpha$ is sufficiently close to $\pi$ then $c$ is (uniformly) quasigeodesic.
\end{cor}

 In higher rank, we do not have an analogue of the inequality \eqref{eq:rank1-qg}, instead, we will be generalizing the corollary. However, {\em angles} in the corollary will be replaced with {\em $\zeta$-angles}.  We will show (in a String of Diamonds Theorem, theorem 
\ref{thm:string}) that if a piecewise-geodesic path $c$ in $X$ has sufficiently long edges and $\zeta$-angles between consecutive segments  
sufficiently close to $\pi$, then $c$ is $M$-Morse for a suitable Morse datum.

\medskip 

In the following,
we consider finite or infinite sequences $(x_n)$ of points in $X$. 

\begin{dfn}[Straight and spaced sequence]
We call a sequence $(x_n)$ 
{\em $(\Theta,\eps)$-straight} 
if the segments $x_nx_{n+1}$ are $\Theta$-regular 
and 
\begin{equation*}
\angle_{x_n}^{\zeta}(x_{n-1},x_{n+1})\geq\pi-\eps
\end{equation*}
for all $n$. 
We call it {\em $l$-spaced} if the segments $x_nx_{n+1}$ 
have length $\geq l$.  
\end{dfn}
Note that every straight sequence 
can be extended to a biinfinite straight sequence. 

Straightness is a local condition.
The goal of this section is to prove the following 
local-to-global result 
asserting that sufficiently straight and spaced sequences 
satisfy a higher rank version of the Morse Lemma 
(for quasigeodesics in hyperbolic space). 

\begin{thm}[Morse Lemma for straight spaced sequences]
\label{thm:locstrimplcoastrseq}
For $\Theta,\Theta',\de$ there exist $l,\eps$ such that:

Every $(\Theta,\eps)$-straight $l$-spaced sequence $(x_n)$
is $\de$-close to a parallel set 
$P(\tau_-,\tau_+)$ with simplices $\tau_{\pm}$ of type $\tau_{mod}$, 
and it moves from $\tau_-$ to $\tau_+$ in the sense that 
its nearest point projection $\bar x_n$ to $P(\tau_-,\tau_+)$ satisfies 
\begin{equation}
\label{eq:movstr}
\bar x_{n\pm m}\in V(\bar x_n,\st_{\Theta'}(\tau_{\pm}))
\end{equation}
for all $n$ and $m\geq1$. 
\end{thm}
\begin{rem}[Global spacing]
\label{rem:globsp}
1. As a corollary of this theorem, we will show that straight spaced sequences are quasigeodesic: 
\begin{equation*}
d(x_n,x_{n+m})\geq clm-2\de
\end{equation*}
with a constant $c=c(\Theta')>0$. See Corollary \ref{cor:retractions}.  In particular, by interpolating the sequence $(x_n)$ via geodesic 
segments we obtain a Morse quasigeodesic in $X$. 

2. Theorem \ref{thm:locstrimplcoastrseq} is a higher-rank generalization of two familiar facts from geometry of Gromov-hyperbolic geodesic metric spaces:  The fact that local quasigeodesics (with suitable parameters) are global quasigeodesics and the Morse lemma stating that 
quasigeodesics stay uniformly close to geodesics. In the higher rank, quasigeodesics, of course, need not be close to geodesics, but, instead (under the straightness assumption), are close to diamonds/Weyl cones/parallel sets. 

3.  One can obviously strengthen the Corollary \ref{cor:rank1} by stating that for each $\eps< \pi$ there exists $L_0(\eps)$ such that if $\al\ge \pi- \eps$ and $L\ge L_0(\eps)$ then  $c$ is a uniform quasigeodesic in $X$. A similar strengthening is false for symmetric spaces of rank $\ge 2$. For instance, when $
W\cong S_3$ and $\eps=2\pi/3$,  then no matter what $\Theta, \Theta'$ and $l$ are, the conclusion of 
Theorem \ref{thm:locstrimplcoastrseq} fails already for sequences contained in a single flat. 
\end{rem}

In order to prove the theorem, 
we start by considering half-infinite sequences 
and prove that they keep moving away from an ideal simplex 
of type $\tau_{mod}$
if they do so initially. 

\begin{dfn}[Moving away from an ideal simplex]
Given a face $\tau\subset \tits X$ of type $\tau_{mod}$ and distinct points $x, y\in X$, define the angle 
$$
\angle_{x}^{\zeta}(\tau, y):= \angle_{x}(z, y)
$$ 
where $z$ is a point (distinct from $x$) on the geodesic ray $x\xi$, where $\xi\in \tau$ is the point of type $\zeta$. 

We say that a sequence $(x_n)$ 
{\em moves $\eps$-away} 
from a simplex $\tau$ of type $\tau_{mod}$ if  
\begin{equation*}
\angle_{x_n}^{\zeta}(\tau,x_{n+1})\geq\pi-\eps
\end{equation*}
for all $n$. 
\end{dfn}

\begin{lem}[Moving away from ideal simplices]
\label{lem:keepsmoving}
For small $\eps$ and large $l$, 
$\eps\leq\eps_0$ and $l\geq l(\eps,\Theta)$, 
the following holds:

If the sequence $(x_n)_{n\geq0}$
is $(\Theta,\eps)$-straight $l$-spaced 
and if 
\begin{equation*}
\angle_{x_0}^{\zeta}(\tau,x_1)\geq\pi-2\eps,
\end{equation*}
then $(x_n)$ moves $\eps$-away from $\tau$. 
\end{lem}
\proof
By Lemma~\ref{lem:distangdec}(ii), 
the unit speed geodesic segment $c:[0,t_1]\to X$ from $p(0)$ to $p(1)$ 
moves $\eps(d(2\eps))$-away from $\tau$ at all times, 
and $\eps'(2\eps,\Theta,l)$-away at times $\geq l$, 
which includes the final time $t_1$. 
For $l(\eps,\Theta)$ sufficiently large, 
we have $\eps'(2\eps,\Theta,l)\leq\eps$. 
Then $c$ moves $\eps$-away from $\tau$ at time $t_1$, 
which means that 
$\angle_{x_1}^{\zeta}(\tau,x_0)\leq\eps$. 
Straightness at $x_1$ and the triangle inequality yield that again
$\angle_{x_1}^{\zeta}(\tau,x_2)\geq\pi-2\eps$. 
One proceeds by induction. 
\qed

\medskip
Note that there do exist simplices $\tau$ 
satisfying the hypothesis of the previous lemma. 
For instance, 
one can extend the initial segment $x_0x_1$ backwards to infinity 
and choose $\tau=\tau(x_1x_0)$. 

\medskip
Now we look at {\em biinfinite} sequences. 

We assume in the following that $(x_n)_{n\in\Z}$
is $(\Theta,\eps)$-straight $l$-spaced
for small $\eps$ and large $l$. 
As a first step, 
we study the asymptotics of such sequences 
and use the argument for Lemma~\ref{lem:keepsmoving} 
to find a pair of opposite ideal simplices $\tau_{\pm}$ 
such that $(x_n)$ moves from $\tau_-$ towards $\tau_+$. 

\begin{lem}[Moving towards ideal simplices]
\label{lem:locstrcpclpar}
For small $\eps$ and large $l$, 
$\eps\leq\eps_0$ and $l\geq l(\eps,\Theta)$, 
the following holds: 

There exists a pair of opposite simplices $\tau_{\pm}$ 
of type $\tau_{mod}$ 
such that the inequality
\begin{equation}
\label{ineq:tauminusplusatnb}
\angle_{x_n}^{\zeta}(\tau_{\mp},x_{n\pm1})\geq\pi-2\eps
\end{equation}
holds for all $n$.
\end{lem}
\proof
1. For every $n$ define a compact set 
$C^\mp_n\subset\Flag(\tau_{mod})$
$$
C^\pm_n=\{\tau_\pm : \angle_{x_n}^{\zeta}(\tau_{\pm},x_{n\mp 1})\geq\pi-2\eps\}. 
$$
As in the proof of Lemma~\ref{lem:keepsmoving}, 
straightness at $x_{n+1}$ implies that 
$C^-_n\subset C^-_{n+1}$. 
Hence the family $\{C^-_n\}_{n\in \Z}$ form a nested sequence of nonempty compact subsets
and therefore have nonempty intersection containing a simplex $\tau_-$.  
Analogously, there exists a simplex $\tau_+$ which belongs to $C^+_n$ for all $n$.  

2. It remains to show that the simplices $\tau_-, \tau_+$ are antipodal. Using straightness and the triangle inequality, we see that 
\begin{equation*}
\label{ineq:tauminusplusatn}
\angle_{x_n}^\zeta(\tau_-,\tau_+)\geq\pi-5\eps 
\end{equation*}
for all $n$. Hence, if $5\eps<\eps(\zeta)$, then the simplices $\tau_-, \tau_+$ are antipodal in view of Remark \ref{rem:antip}. \qed

\medskip
The pair of opposite simplices $(\tau_-, \tau_+)$ which we found determines a parallel set in $X$. 
The second step is to show that $(x_n)$ is uniformly close to it. 

\begin{lem}[Close to parallel set]
\label{lem:closeparsetb}
For small $\eps$ and large $l$, 
$\eps\leq\eps(\de)$ and $l\geq l(\Theta,\de)$, 
the sequence $(x_n)$ is $\de$-close to $P(\tau_-,\tau_+)$. 
\end{lem}
\proof The statement follows from the combination of 
the inequality \eqref{ineq:tauminusplusatn} (in the second part of the proof of Lemma \ref{lem:locstrcpclpar}) and Lemma~\ref{lem:distangcontr}.
\qed

\medskip
The third and final step of the proof 
is to show that the nearest point projection $(\bar x_n)$ of $(x_n)$ 
to $P(\tau_-,\tau_+)$ 
moves from $\tau_-$ towards $\tau_+$.

\begin{lem}[Projection moves towards ideal simplices]
\label{lem:projpathstrb}
For small $\eps$ and large $l$, 
$\eps\leq\eps_0$ and $l\geq l(\eps,\Theta,\Theta')$, 
the segments $\bar x_n\bar x_{n+1}$ are $\Theta'$-regular
and 
\begin{equation*}
\label{ineq:tauminusplusatnproj}
\angle_{\bar x_n}^{\zeta}(\tau_-,\bar x_{n+1})=\pi
\end{equation*}
for all $n$. 
\end{lem}
\proof
By the previous lemma, 
$(x_n)$ is $\delta_0$-close to $P(\tau_-,\tau_+)$
if $\eps_0$ is sufficiently small and $l$ is sufficiently large. 
Since $x_nx_{n+1}$ is $\Theta$-regular, 
the triangle inequality for $\De$-lengths yields that 
the segment $\bar x_n\bar x_{n+1}$ is $\Theta'$-regular, 
again if $l$ is sufficiently large.

Let $\xi_+$ denote the ideal endpoint of the ray extending this segment, 
i.e.\ $\bar x_{n+1}\in\bar x_n\xi_+$. 
Then $x_{n+1}$ is $2\delta_0$-close to the ray $x_n\xi_+$. 
We obtain that 
\begin{equation*}
\tangle^{\zeta}(\tau_-,\xi_+)\geq
\angle^{\zeta}_{x_n}(\tau_-,\xi_+)\simeq
\angle^{\zeta}_{x_n}(\tau_-,x_{n+1})\simeq\pi
\end{equation*}
where the last step follows from inequality (\ref{ineq:tauminusplusatnb}). 
The discreteness of Tits distances between ideal points 
of fixed type $\zeta$ implies that in fact
\begin{equation*}
\tangle^{\zeta}(\tau_-,\xi_+)=\pi, 
\end{equation*}
i.e.\ the ideal points $\zeta(\tau_-)$ and $\zeta(\xi_+)$ 
are antipodal. 
But the only simplex opposite to $\tau_-$ in 
$\geo P(\tau_-,\tau_+)$ is $\tau_+$, 
so $\tau(\xi_+)=\tau_+$ and 
\begin{equation*}
\angle_{\bar x_n}^{\zeta}(\tau_-,\bar x_{n+1})=
\angle_{\bar x_n}^{\zeta}(\tau_-,\xi_+)
=\pi,
\end{equation*}
as claimed.
\qed

\medskip
{\em Proof of Theorem~\ref{thm:locstrimplcoastrseq}.}
It suffices to consider biinfinite sequences. 

The conclusion of 
Lemma~\ref{lem:projpathstrb} 
is equivalent to 
$\bar x_{n+1}\in V(\bar x_n,\st_{\Theta'}(\tau_+))$. 
Combining Lemmas~\ref{lem:closeparsetb} and~\ref{lem:projpathstrb}, 
we thus obtain the theorem for $m=1$. 

The convexity of $\Theta'$-cones, cf.\ Proposition~\ref{prop:thconeconv}, 
implies that 
\begin{equation*}
V(\bar x_{n+1},\st_{\Theta'}(\tau_+))\subset 
V(\bar x_n,\st_{\Theta'}(\tau_+)),
\end{equation*}
and the assertion follows for all $m\geq1$ by induction.  
\qed

\begin{rem}
\label{rem:xinotconv}
The conclusion of the theorem implies flag-convergence  
$x_{\pm n}\to\tau_{\pm}$ as $n\to+\infty$. 
However, 
the sequences $(x_n)_{n\in \pm \N}$ do in general not converge at infinity,
but accumulate at compact subsets of $\st_{\Theta'}(\tau_{\pm})$. 
\end{rem}

\subsection{Lipschitz retractions to straight  paths}

Consider a (possibly infinite) closed interval $J$ in $\R$; we will assume that $J$ has integer or infinite bounds.  
Suppose that $p: J\cap \Z\to P=P(\tau_-, \tau_+)\subset X$ is an $l$-separated, $\lambda$-Lipschitz, $(\Theta,0)$-straight coarse sequence 
pointing away from $\tau_-$ and towards $\tau_+$. We extend $p$ to a piecewise-geodesic map $p: J\to P$ by sending intervals $[n, n+1]$ to geodesic segments $p(n)p(n+1)$ via affine maps. We retain the name $p$ for the extension. 

\begin{lem}
There exists $L=L(l,\lambda,\Theta)$ and an  $L$-Lipschitz retraction of $X$ to $p$, i.e., an $L$-Lipschitz map $r: X\to J$ so that $r\circ p=Id$. 
In particular, $p: J\cap \Z\to X$ is a $(\bar L, \bar A)$-quasigeodesic, where $\bar L, \bar A$ depend only on $l, \lambda, \Theta$.  
\end{lem}
\proof It suffices to prove existence of a retraction. Since $P$ is convex in $X$, 
it suffices to construct a map $P\to J$. Pick a generic point $\xi=\xi_+\in \tau_+$ and let $b_\xi: P\to \R$  denote the Busemann function normalized so that $b_\xi(p(z))=0$ for some $z\in J\cap \Z$.  Then the $\Theta$-regularity assumption on $p$ implies that the slope of the piecewise-linear function $b_\xi\circ p: J\to \R$ is strictly positive, bounded away from $0$. The assumption that $p$ is $l$-separated $\la$-Lipschitz implies that
$$
l\le |p'(t)|\le \la  
$$
for each $t$ (where the derivative exists). The straightness assumption on $p$ implies that the function $h:= b_\xi\circ p: J\to \R$ is strictly increasing. By combining these observations, we conclude that $h$ is an $L$-biLipschitz homeomorphism for some 
$L=L(l,\la,\Theta)$. Lastly, we define 
$$
r: P\to J, \quad r=h^{-1}\circ b_\xi. 
$$ 
Since $b_\xi$ is $1$-Lipschitz, the map $r$ is $L$-Lipschitz. By the construction, $r\circ p=Id$. \qed 

\begin{cor}\label{cor:retractions}
Suppose that $p: J\cap \Z\to X$ is a $l$-spaced, $\lambda$-Lipschitz, $(\Theta,\eps)$-straight  sequence. 
Pick some $\Theta'$ such that $\Theta\subset int(\Theta')$ and let $\delta=\delta(l, \Theta, \Theta', \eps)$ 
be the constant as in Theorem \ref{thm:locstrimplcoastrseq}. Then for 
$L=L(l-2\delta,\la+2\delta,\Theta')$  we have: 

1. There exists an $(L,2\delta)$-coarse Lipschitz retraction $X\to J$. 

2. The map $p$ is a $(\Theta', D', L',A')$-quasigeodesic with $D', L', A'$ depending only on  $l, \la, \Theta, \Theta', \eps$.  
\end{cor}
\proof The statement immediately follows the above lemma combined with Theorem \ref{thm:locstrimplcoastrseq}.  \qed 

\medskip
Reformulating in terms of piecewise-geodesic paths, we obtain 

\begin{thm}
[String of diamonds theorem]  \label{thm:string}
For any pair of Weyl convex subsets $\Theta< \Theta'$ and a number $D\ge 0$ 
there exist positive numbers $\eps$, $S$, $L$, $A$ depending on the datum 
$(\Theta, \Theta',D)$ such that the following holds. 

Suppose that $c$ is an arc-length parameterized piecewise-geodesic path (finite or infinite) 
in $X$ obtained by concatenating geodesic segments $x_i x_{i+1}$ satisfying for all $i$:

1. Each segment $x_i x_{i+1}$ is $\Theta$-regular and has length $\ge S$. 

2. 
$$
\zangle_{x_i}(x_{i-1}, x_{i+1})\ge \pi -\eps. 
$$
Then the path $c$ is $(\Theta', D, L, A)$-Morse.  
\end{thm}

\subsection{Local Morse quasigeodesics}
\label{sec:morseqg}

According to Theorem~\ref{thm:string}, 
sufficiently straight and spaced straight piecewise-geodesic paths are Morse.  
In this section we will now prove that, conversely, 
the Morse property implies straightness in a suitable sense, 
namely that for sufficiently spaced quadruples 
the associated midpoint triples are arbitrarily straight. 
(For the quadruples themselves this is in general not true.)

\begin{dfn}[Quadruple condition]
\label{dfn:quad}
For points $x, y\in X$ we let $\midp(x,y)$ denote the midpoint of the geodesic segment $xy$. 
A map $p:I\to X$ satisfies the 
{\em $(\Theta,\eps,l,s)$-quadruple condition}
if for all $t_1,t_2,t_3,t_4\in I$
with $t_2-t_1,t_3-t_2,t_4-t_3\geq s$ 
the triple of midpoints
\begin{equation*}
(\midp(t_1,t_2),\midp(t_2,t_3),\midp(t_3,t_4))
\end{equation*}
is $(\Theta,\eps)$-straight and $l$-spaced. 
\end{dfn}
\begin{prop}[Morse implies quadruple condition]
\label{prop:morseimplquad}
For $L,A,\Theta,\Theta',D,\eps,l$ 
exists a scale $s=s(L,A,\Theta,\Theta',D,\eps,l)$ 
such that every $(\Theta,D,L,A)$-Morse quasigeodesic 
satisfies the $(\Theta',\eps,l,s')$-quadruple condition for every $s'\ge s$.  
\end{prop}
\proof
Let $p:I\to X$ be an $(L,A,\Theta,D)$-Morse quasigeodesic, 
and let $t_1,\dots,t_4\in I$
such that $t_2-t_1,t_3-t_2,t_4-t_3\geq s$. 
We abbreviate $p_i:=p(t_i)$ and $m_i=\midp(p_i,p_{i+1})$. 

Regarding straightness, 
it suffices to show that 
the segment $m_2m_1$ is $\Theta'$-regular 
and that $\angle^{\zeta}_{m_2}(p_2,m_1)\leq\frac{\eps}{2}$
provided that $s$ is sufficiently large 
in terms of the given data. 

By the Morse property,
there exists a diamond
$\diamo_{\Theta}(x_1,x_3)$
such that $d(x_1,p_1),d(x_3,p_3)\leq D$ 
and $p_2\in N_D(\diamo_{\Theta}(x_1,x_3))$.
The diamond spans a unique parallel set $P(\tau_-,\tau_+)$. 
(Necessarily, 
$x_3\in V(x_1,\st_{\Theta}(\tau_+))$ and 
$x_1\in V(x_3,\st_{\Theta}(\tau_-))$.)

We denote by $\bar p_i$ and $\bar m_i$
the projections of $p_i$ and $m_i$ to the parallel set. 

We first observe 
that $m_2$ (and $m_3$) is arbitrarily close to the parallel set 
if $s$ is large enough. 
If this were not true,
a limiting argument would produce a geodesic line 
at strictly positive finite Hausdorff distance 
$\in(0,D]$ 
from $P(\tau_-,\tau_+)$ 
and asymptotic to ideal points in 
$\st_{\Theta}(\tau_{\pm})$. 
However,
all lines asymptotic to ideal points in
$\st_{\Theta}(\tau_{\pm})$ 
are contained in $P(\tau_-,\tau_+)$. 

Next,
we look at the directions of the segments $\bar m_2\bar m_1$ 
and $\bar m_2\bar p_2$
and show that they have the same $\tau$-direction. 
Since $\bar p_2$ is $2D$-close to 
$V(\bar p_1,\st_{\Theta}(\tau_+))$, 
we have that the point 
$\bar p_1$ is $2D$-close to 
$V(\bar p_2,\st_{\Theta}(\tau_-))$, 
and hence also 
$\bar m_1$ is $2D$-close to 
$V(\bar p_2,\st_{\Theta}(\tau_-))$.
Therefore, 
$\bar p_1,\bar m_1\in V(\bar p_2,\st_{\Theta'}(\tau_-))$ 
if $s$ is large enough. 
Similarly, 
$\bar m_2\in V(\bar p_2,\st_{\Theta'}(\tau_+))$
and hence 
$\bar p_2\in V(\bar m_2,\st_{\Theta'}(\tau_-))$. 
The convexity of $\Theta'$-cones, 
see Proposition~\ref{prop:thconeconv}, 
implies that also 
$\bar m_1\in V(\bar m_2,\st_{\Theta'}(\tau_-))$. 
In particular,
$\angle^{\zeta}_{\bar m_2}(\bar p_2,\bar m_1)=0$ 
if $s$ is sufficiently large. 

Since $m_2$ is arbitrarily close to the parallel set 
if $s$ is sufficiently large, 
it follows by another limiting argument that 
$\angle^{\zeta}_{m_2}(p_2,m_1)\leq\frac{\eps}{2}$
if $s$ is sufficiently large. 

Regarding the spacing, 
we use that 
$\bar m_1\in V(\bar p_2,\st_{\Theta'}(\tau_-))$
and 
$\bar m_2\in V(\bar p_2,\st_{\Theta'}(\tau_+))$.
It follows that 
\begin{equation*}
d(\bar m_1,\bar m_2) \geq c\cdot (d(\bar m_1,\bar p_2)+d(\bar p_2,\bar m_2))
\end{equation*}
with a constant $c=c(\Theta')>0$, 
and hence that 
$d(m_1,m_2)\geq l$ 
if $s$ is sufficiently large. 
\qed

\medskip
Theorem~\ref{thm:locstrimplcoastrseq} and Proposition~\ref{prop:morseimplquad} 
tell that the Morse property for quasigeodesics is equivalent 
to straightness (of associated spaced sequences of points). 
Since straightness is a local condition, 
this leads to a local to global result for Morse quasigeodesics, 
namely that the Morse property holds globally 
if it holds locally up to a sufficiently large scale. 

\begin{dfn}[Local Morse quasigeodesic]
\label{defn:locmqg} 
An $S$-local $(\Theta,D,L,A)$-{\em Morse quasigeo\-de\-sic} in $X$ 
is a map $p:I\to X$ 
such that for all $t_0$ the subpath 
$p|_{[t_0,t_0+S]}$ is a $(\Theta,D,L,A)$-Morse quasigeodesic. 
\end{dfn}
Note that local Morse quasigeodesics are uniformly coarse Lipschitz.

\begin{thm}[Local-to-global principle for Morse quasigeodesics]
\label{thm:locglobmqg}
For $L,A,\Theta,\Theta',D$ exist $S,L',A',D'$ such that 
every $S$-local 
$(\Theta,D,L,A)$-local Morse quasigeo\-de\-sic in $X$ 
is an $(\Theta',D',L',A')$-Morse quasigeodesic.  
\end{thm}
\proof
We choose an auxiliary Weyl convex subset $\Theta''$ 
such that $\Theta< \Theta''< \Theta'$.  

Let $p:I\to X$ be an $S$-local $(\Theta,D,L,A)$-local Morse quasigeo\-de\-sic. 
We consider its coarsification on a (large) scale $s$ 
and the associated midpoint sequence, 
i.e.\ we put $p^s_n=p(ns)$ and $m^s_n=\midp(p^s_n,p^s_{n+1})$. 
Whereas the coarsification itself 
does in general not become arbitrarily straight 
as the scale $s$ increases, 
this is true for its midpoint sequence
due to Proposition~\ref{prop:morseimplquad}. 
We want it to be sufficiently straight and spaced 
so that we can apply to it the Morse Lemma from 
Theorem~\ref{thm:locstrimplcoastrseq}. 
Therefore we first fix an auxiliary constant $\de$, 
and further auxiliary constants $l,\eps$ 
as determined by Theorem~\ref{thm:locstrimplcoastrseq} 
in terms of $\Theta',\Theta''$ and $\de$. 
Then 
Proposition~\ref{prop:morseimplquad} applied to the 
$(\Theta,D,L,A)$-Morse quasigeodesics $p|_{[t_0,t_0+S]}$ 
yields that 
$(m^s_n)$ is $(\Theta'',\eps)$-straight and $l$-spaced 
if $S\geq 3s$ and the scale $s$ is large enough 
depending on $L,A,\Theta,\Theta'',D,\eps,l$. 

Now we can apply Theorem~\ref{thm:locstrimplcoastrseq} to $(m^s_n)$. 
It yields 
a nearby sequence $(\bar m^s_n)$, 
$d(\bar m^s_n,m^s_n)\leq\de$, 
with the following property:
For all $n_1<n_2<n_3$ 
the segments $\bar m^s_{n_1}\bar m^s_{n_3}$ 
are uniformly regular 
and the points $m^s_{n_2}$ are $\de$-close to the diamonds
$\diamo_{\Theta'}(\bar m^s_{n_1},\bar m^s_{n_3})$. 

Since the subpaths $p|_{[ns,(n+1)s]}$ filling in $(p^s_n)$ 
are $(L,A)$-quasigeodesics (because $S\geq s$), 
and it follows that 
for all $t_1,t_2\in I$ 
the subpaths $p|_{[t_1,t_2]}$ are 
$D'$-close to $\Theta'$-diamonds 
with $D'$ depending on $L,A,s$. 

The conclusion of Theorem~\ref{thm:locstrimplcoastrseq} 
also implies a global spacing for the sequence $(m^s_n)$, 
compare Remark~\ref{rem:globsp}, i.e.\ 
$d(m^s_n,m^s_{n'})\geq c\cdot |n-n'|$
with a positive constant $c$ depending on $\Theta',l$. 
Hence $p$ is a global $(L',A')$-quasigeodesic 
with $L',A'$ depending on 
$L,A,s,c$. 

Combining this information, 
we obtain that $p$ is an $(\Theta',D',L',A')$-Morse quasigeodesic 
for certain constants $L',A'$ and $D'$ 
depending on $L,A,\Theta,\Theta'$ and $D$, 
provided that the scale $S$ is sufficiently large 
in terms of the same data. 
\qed

\subsection{Local-to-global principle for Morse maps}\label{sec:locglobmm}

We now deduce from our local-to-global result for Morse quasigeodesics, 
 Theorem~\ref{thm:locglobmqg}, 
a local-to-global result for Morse embeddings. 

\medskip
We restrict to the setting of maps of Gromov-hyperbolic $(l,a)$-quasigeodesic metric spaces $Z$ to symmetric spaces $X$.


\begin{dfn}[Local Morse embedding]
We call a map $f:Z\to X$ an $S$-local $(\Theta,D,L,A)$-Morse map 
if   for any $(l,a)$-quasigeodesic $q: I \to Z$ 
defined on an interval $I$ of length $\leq S$ 
the image path $f\circ q$ is a  
$(\Theta,D,L,A)$-Morse quasigeodesic in $X$. 
\end{dfn}

\begin{thm}[Local-to-global principle for Morse embeddings of Gromov hyperbolic spaces]
\label{thm:locglobmqiembhyp0}
For $l,a,L,A,\Theta,\Theta',D$ exists a scale $S$ and a datum $(D',L',A')$ 
such that every $S$-local 
$(\Theta,D,L,A)$-Morse embedding 
from an $(l,a)$-quasigeodesic  Gromov hyperbolic space into $X$ 
is a $(\Theta',D',L',A')$-Morse embedding. 
\end{thm}
\proof 
Let $f:Z\to X$ denote the local Morse embedding. 
It sends every $(l,a)$-quasigeodesic $q:I\to Z$ 
to an $S$-local $(\Theta,D,L,A)$-Morse quasigeodesic $p=f\circ q$ in $X$. 
By Theorem~\ref{thm:locglobmqg}, 
$p$ is $(L',A',\Theta',D')$-Morse 
if $S\ge S(l,a,L,A,\Theta,\Theta',D)$, 
where $L',A',D'$ depend on the given data. 
\qed

Below is a reformulation of this theorem in the case of geodesic Gromov-hyperbolic spaces. 

Let $Z$ be a $\delta$-hyperbolic geodesic space. An $R$-ball $B(z,R)$ in $Z$ need not be convex, but it is $\delta$-quasiconvex. 
In particular, the restriction of the metric from $Z$ to $B(z,R)$ results in a $(1,\delta)$-quasigeodesic metric space.

\begin{thm}[Local-to-global principle for Morse embeddings of geodesic spaces]
\label{thm:locglobmqiembhyp}
For $L,A,\Theta,\Theta',D,\delta$ exists a scale $R$ and a datum $(D',L',A')$ 
such that if $Z$ is a $\delta$-hyperbolic geodesic metric space and the restriction of $f$ to any $R$-ball is $(\Theta,D,L,A,1,\delta)$-Morse, 
then $f: Z\to X$ is  $(\Theta',D',L',A')$-Morse. 
\end{thm}

\section{Group-theoretic applications} \label{sec:group-applications}

As a consequence of the local-to-global criterion for Morse maps, in this section we establish that the Morse property for isometric group actions is an open condition. Furthermore, for two nearby Morse actions, the actions on their $\tau_{mod}$-limit sets are also close, i.e.\ conjugate by an equivariant homeomorphism close to identity. In view of the equivalence of Morse property with the asymptotic properties discussed earlier,  this implies structural stability for asymptotically embedded groups.
Another corollary of the local-to-global result is the algorithmic recognizability of Morse actions.

We conclude the section by illustrating our technique by constructing Morse-Schottky actions of free groups on higher rank symmetric spaces.

\subsection{Stability of Morse actions}

We consider isometric actions $\Ga\acts X$ 
of finitely generated groups. 
\begin{dfn}[Morse action]
\label{dfn:morseact}
We call an action $\Ga\acts X$ {\em $\Theta$-Morse} 
if one (any) orbit map $\Ga\to\Ga x\subset X$ 
is a $\Theta$-Morse embedding 
with respect to a(ny) word metric on $\Ga$. 
We call an action $\Ga\acts X$ {\em $\tau_{mod}$-Morse} if it is $\Theta$-Morse 
for some $\tau_{mod}$-Weyl convex compact subset $\Theta\subset\ost(\tau_{mod})$. 
\end{dfn}

\begin{rem}[Morse actions are $\taumod$-regular and undistorted]
\label{rem:morseacundist}
(i) It follows immediately from the definition of Morse quasigeodesics that $\Theta$-Morse actions 
are $\tau_{mod}$-regular for the simplex type $\tau_{mod}$ determined by $\Theta$. 

(ii) Morse subgroups of $G$ are {\em undistorted} in the sense 
that the orbit maps are quasi-isometric embeddings. 
In \cite{bordif} we prove that Morse subgroups of $G$ satisfy a stronger property: They are {\em coarse Lipschitz retracts} of $G$.  This retraction property is stronger than nondistortion: Every finitely generated 
subgroup which is a coarse retract of $G$ is undistorted in $G$, but there are  examples of undistorted subgroups which are not coarse  retracts. 
For instance, the group $\Phi:= F_2\times F_2$ admits an undistorted embedding in the isometry group of $X=\H^2\times \H^2$. On the other hand, pick  an epimorphism 
$\phi: F_2\to \Z$   and define the subgroup $\Ga< \Phi$ as the kernel of the homomorphism 
$$
(\gamma_1, \gamma_2) \mapsto \phi(\gamma_1) - \phi(\gamma_2). 
$$
Then $\Gamma$ is a finitely generated undistorted subgroup of $\Phi$ (see e.g. \cite[Theorem 2]{OS}), but is not finitely presented (see e.g. \cite{BR}). 
Hence, $\Ga< G = \Isom(\H^2)\times \Isom(\H^2)$ is undistorted but is not a  coarse Lipschitz retract. 
\end{rem}

We denote by $\Hom_{\taumod}(\Ga,G)\subset\Hom(\Ga,G)$ 
the subset of $\taumod$-Morse actions $\Ga\acts X$.

By analogy with {\em local Morse quasigeodesics}, we define {\em local Morse group actions} $\rho: \Ga\acts X$ of a hyperbolic group 
(with a fixed finite generating set):  

\begin{dfn}
An action $\rho$ is called  $S$-locally {\em $(\Theta, D, L, A)$-locally Morse}, or 
{\em $(\Theta, D, L,A)$-locally Morse on the scale $S$}, with respect to a base-point $x\in X$, if the 
orbit map $\Ga\to \Ga\cdot x\subset X$ induces an $S$-local  $(\Theta, D, L, A)$-local Morse embedding of the Cayley graph of $\Ga$. 
\end{dfn}

According to our local-to-global result for Morse embeddings, 
see Theorem~\ref{thm:locglobmqiembhyp}, 
an action of a word hyperbolic group is Morse 
if and only if it is local Morse on a sufficiently large scale. 
Since this is a finite condition, 
it follows that the Morse property is stable under perturbation of the action:

\begin{thm}[Morse is open for word hyperbolic groups]
\label{thm:morsestab}
For any word hyperbolic group $\Ga$ 
the subset $\Hom_{\taumod}(\Ga,G)$ is open in $\Hom(\Ga,G)$. More precisely, if $\rho\in \Hom_{\taumod}(\Ga,G)$ is $M$-Morse with respect to a base-point $x\in X$ then there exists a neighborhood of $\rho$ in   $\Hom(\Ga,G)$ consisting entirely of $M'$-Morse representations with respect to $x$, where $M'$ depends only on $M$. 
\end{thm}
\proof
Let $\rho:\Ga\acts X$ be a Morse action. 
We fix a word metric on $\Ga$ and a base point $x\in X$. 
Then there exist data $M=(L, A, \Theta, D)$ 
such that the orbit map $\Ga\to\Ga x\subset X$ extends to  
a $(\Theta, D, L, A)$-Morse map of the Cayley graph $Y$ on $\Ga$.

We relax the Morse parameters slightly, 
i.e.\ we consider $(L,A,\Theta,D)$-Morse quasigeodesics 
as $(L,A+1,\Theta,D+1)$-Morse quasigeodesics
satisfying strict inequalities. 
For every scale $S$, the orbit map $\Ga \to \Ga x\subset X$, defines 
an $(L,A+1,\Theta,D+1,S)$-{\em local} Morse embedding $Y\to X$. 
Due to $\Ga$-equivariance, 
this is a finite condition in the sense 
that it is equivalent to a condition involving 
only finitely many orbit points. 
Since we relaxed the Morse parameters, 
the same condition is satisfied by all actions sufficiently close to $\rho$.

Theorem~\ref{thm:locglobmqiembhyp} provides a scale $S$ 
such that all $S$-local $(\Theta,D+1, L, A+1)$-Morse embeddings $Y\to X$ 
are $M'$-Morse for some Morse datum $M'$ depending only on $(L,A+1,\Theta,D+1,S)$. 
It follows that actions sufficiently close to $\rho$ are $\taumod$-Morse. 
\qed

\begin{cor}
For every hyperbolic group $\Ga$ the space of {\em faithful} Morse representations 
$$\Hom_{inj,\taumod}(\Ga,G)$$ is open in $\Hom_{\taumod}(\Ga,G)$. 
\end{cor}
\proof Every hyperbolic group $\Ga$ has the unique maximal finite normal subgroup 
$\Phi\triangleleft \Ga$ (if $\Ga$ is nonelementary then $\Phi$ is the kernel of the action of $\Ga$ on $\geo \Ga$). Since Morse actions are properly discontinuous, the kernel of every Morse representation $\Ga\to G$ is contained in $\Phi$. Since $\Hom(\Phi, G)/G$ is finite, it follows that the set of faithful Morse representations is open in $\Hom_{\taumod}(\Ga,G)$. \qed

\medskip
The result on the openness of the Morse condition 
for actions of word hyperbolic groups, 
cf.\ Theorem~\ref{thm:morsestab}, 
can be strengthened in the sense that 
the asymptotics of Morse actions vary continuously:

\begin{thm}[Morse actions are structurally stable]
\label{thm:strcstab}
The boundary map at infinity of a Morse action 
depends continuously on the action. 
\end{thm}
\proof According to Theorem~\ref{thm:morsestab} 
nearby actions are uniformly Morse. 
The assertion therefore follows from the fact 
that the ends of Morse quasirays
vary uniformly continuously, 
cf.\ Lemma~\ref{lem:cont-dependence}.
\qed

\begin{rem}
(i)
Note that since the boundary maps at infinity are embeddings, 
the $\Ga$-actions on the $\tau_{mod}$-limit sets 
are topologically conjugate to each other and, for nearby actions, 
by a homeomorphism close to the identity.

(ii)
In rank one, our argument yields a different proof for 
Sullivan's Structural Stability Theorem \cite{Sullivan}
for convex cocompact group actions on rank one symmetric spaces. Other proofs can be found in \cite{Labourie, GW} 
(for Anosov subgroups in higher rank), \cite{Corlette, Izeki, Bowditch-stab} for rank one symmetric spaces.  
\end{rem}

\medskip
Our next goal is to extend the topological conjugation from the limit set to the domains of proper discontinuity. Recall that in \cite{coco15} we constructed domains of proper discontinuity and cocompactness for $\taumod$-Morse group actions on flag-manifolds $\Flagn=G/P_{\numod}$. Such domains depend on a certain auxiliary datum, a {\em balanced thickening} 
$\Th\subset W$, which is a $W_{\taumod}$-left invariant subset satisfying certain conditions; see \cite[sect. 3.4]{coco15}. 
Let $\numod\subset \simod$ be an $\iota$-invariant face such that $\Th$ is invariant under the action of $W_{\numod}$ via the {\em right} multiplication (this is automatic if 
$\numod=\simod$ since $W_{\simod}=\{e\}$). The thickening $\Th\subset W$ defines a thickening $\Th(\Lat(\Ga))\subset \Flagn$.  
One of the main results of \cite{coco15} (Theorem 1.7) is that each $\taumod$-Morse subgroup $\Ga< G$ acts properly discontinuously and cocompactly on 
 $$
 \Om_{\Th}(\Ga):=\Flagn - \Th(\Lat(\Ga)). 
 $$

\begin{thm}
[Stability of Morse quotient spaces] \label{thm:Omstability}
Suppose that $\rho_n: \Ga\to \rho_n(\Ga)=
\Ga_n< G$ is a sequence of faithful $\taumod$-Morse representations converging to a $\taumod$-Morse embedding $\rho: \Ga\embed G$. Then:

1. The sequence of thickenings $\Th(\Lat(\Ga_n))$ Hausdorff-converges to $\Th(\Lat(\Ga))$. 

2. If $\ga_n\in \Ga$ is a divergent sequence, then, after extraction, the sequence $(\rho_n(\ga_n))$ flag-converges to a point in $\Lat(\Ga)$. 

3. There is a sequence of equivariant diffeomorphisms $h_n: \Om_{\Th}(\Ga)\to \Om_{\Th}(\Ga_n)$ converging to the identity map uniformly on compacts. 

4. In particular, the quotient-orbifolds $\Om_{\Th}(\Ga_n)/\Ga_n$ are diffeomorphic to 
$\Om_{\Th}(\Ga)/\Ga$ for all sufficiently large $n$. 
\end{thm}
\proof 1. First of all, suppose that a sequence $\tau_n\in \Flagt$ converges to $\tau\in \Flagt$. Then, since $\Flagn= G/P_{\numod}$, there is a sequence $g_n\in G$, $g_n\to e$, such that $g_n(\tau)=\tau_n$. Since
$$
g_n(\Th(\tau))=\Th(g_n \tau)= \Th(\tau_n),
$$
it follows that we have Hausdorff-convergence of subsets $\Th(\tau_n)\to \Th(\tau)$. 
Moreover, this convergence of subsets is uniform: There exists $n_0=n(\delta)$ such that if $d(\tau_n,\tau)<\delta$ for all $n\ge n_0$ 
then $d(\Th(\tau_n), \Th(\tau))< \eps=\eps(\delta)$ for all $n\ge n_0$. Here $\eps\to 0$ as $\delta\to 0$. 
Since the sequence of limit sets $\Lat(\Ga_n)$ Hausdorff-converges to $\Lat(\Ga)$, 
it follows that the sequence of thickenings $\Th(\Lat(\Ga_n))$  Hausdorff-converges to $\Th(\Lat(\Ga))$. This proves (1). 

\medskip
2. Consider a sequence of geodesic rays $e \xi_n$ in the Cayley graph $Y$ of $\Ga$ such that $\ga_n$ lies in an $R$-neighborhood of $e\xi_n$ for all $n$. Then, in view of the uniform $M'$-Morse property for the representations $\rho_n$, each point $\rho_n(\ga_n)(x)$ belongs to the $D'$-neighborhood of the Weyl cone $V(x, \st(\tau_n))$, where $\tau_n= \al_n(\xi_n)$, $\al_n: \geo \Ga\to \Lat(\Ga_n)$ is the asymptotic embedding. Thus, by the definition of flag-convergence, the sequences $(\rho_n(\ga_n))$ and $(\tau_n)$ have the same flag-limit in $\Flagt$. By Part 1, the sequence $(\tau_n)$ subconverges to a point in $\Lat(\Ga)$. Hence, the same holds for $(\rho_n(\ga_n))$.

\medskip 
3. The proof of this part is mostly standard, see \cite{Izeki} in the case  when $X$ is a hyperbolic space. 
The quotient orbifold $O=\Om_{\Th}(\Ga)/\Ga$ has a natural  $(\F, G)$-structure where $\F=\Flagn$. 
The orbifold $O$ has finitely many components, let $Z$ be one of them and let $\hat{Z}\subset \Om_{\Th}(\Ga)$ be a component projecting to $Z$.  It suffices to construct maps $h_n$ on each component $\hat{Z}$ and then extend these maps to maps $h_n$ of $\Om_{\Th}(\Ga)$ by $\rho_n$-equivariance.

The covering map $\hat{Z}\to Z$ induces an epimorphism $\phi: \pi_1(Z)\to \Ga_Z$, where $\Ga_Z$ is the $\Ga$-stabilizer of $\hat{Z}$. Let $dev: \tilde{Z}\to  
\hat{Z}\subset \Om_{\Th}(\Ga)$ be the developing map, where $\tilde{Z}\to Z$ is the universal covering. 
By Ehresmann-Thurston holonomy theorem (see \cite{Lok}, \cite{CEG}, \cite{Goldman}, \cite[sect. 7.1]{Kapovich-book}), for all sufficiently large $n$, the homomorphism 
$\phi_n:=\rho_n\circ \phi$ is the holonomy of an  $(\F, G)$-structure on $Z$. Moreover, the developing maps $dev_n: \tilde{Z}\to \F$ converge to $dev$ uniformly on compacts in the $C^\infty$-topology. 
Since $\pi_1(\hat{Z})$ is contained in the kernel of $\phi$, it is also in the kernel of $\phi_n$. Hence, the maps $dev_n$ 
descend to maps $\widehat{dev}_n: \hat{Z}\to \F$. The sequence $\widehat{dev}_n$ still converges to 
the identity embedding $\hat{Z}\embed \F$ uniformly on compacts. Pick a compact fundamental set 
$C\subset \hat{Z}$ for the $\Ga_Z$-action, i.e. a compact subset whose $\Ga$-orbit equals $\hat{Z}$. In view of Part 1 of the theorem,  
$\widehat{dev}_n(C)\subset \Om_{\Th}(\Ga_n)$ for all sufficiently large $n$. Therefore,  
 we can assume that $\widehat{dev}_n(\hat{Z})$ is contained in a component $\hat{Z}_n$ 
of $\Om_{\Th}(\Ga_n)$. By the compactness of the quotient-orbifolds,  $\widehat{dev}_n$ projects to a finite-to-one (smooth) orbi-covering map 
$c_n: Z\to Z_n:= \hat{Z}_n/\rho_n(\Ga_Z)$. Hence, $\widehat{dev}_n: \hat{Z}\to \hat{Z}_n$ is a covering map as well. If $\hat{Z}_n$ were simply-connected, it would follow that 
$\widehat{dev}_n$ is a diffeomorphism as required (and this is how Izeki concludes his proof in \cite{Izeki}). We will prove that 
$\widehat{dev}_n$ is a diffeomorphism by a direct argument. 

Suppose that each $\widehat{dev}_n$ is not injective. Then, by the equivariance of these maps, after extraction, there exist convergent sequences $z_n\to z, z'_n\to z'$ in $\hat{Z}$  and a sequence $\ga_n\in \Ga$ such that 
$$
\rho_n(\ga_n) \widehat{dev}_n(z_n)= \widehat{dev}_n(z'_n), \quad \ga_n(z_n)\ne z'_n.  
$$
If the sequence $(\ga_n)$ were contained in a finite subset of $\Ga$ we would obtain a contradiction with the uniform convergence on compacts $\widehat{dev}_n\to id$ on $\hat{Z}$. Hence, 
after extraction, we may assume that $(\ga_n)$ is a divergent sequence. We, therefore, obtain a dynamical relation between the points $z, z'$ via the sequence $(\rho_n(\ga_n))$. 
According to Part 2, the sequence $(\rho_n(\ga_n))$ flag-accumulates to $\Lat(\Ga)$. The dynamical relation then contradicts fatness of the balanced thickening $\Th$, see 
\cite[sect. 5.2]{coco15} and the proof of Theorem 6.8 in \cite{coco15}. 

We conclude that the maps 
$$
\widehat{dev}_n: \hat{Z}\to \hat{Z}_n
$$
 are diffeomorphisms for all sufficiently large $n$.  Since $\rho_n: \Ga\to \Ga_n$ are isomorphisms, equivariance of the developing maps implies that the maps 
 $h_n: \Om_{\Th}(\Ga)\to \Om_{\Th}(\Ga_n)$  are diffeomorphisms for sufficiently large $n$. 

4. This part is an immediate corollary of Part 3. \qed

\begin{rem}
(i) In the case  when $X$ is a hyperbolic space, the equivariant diffeomorphism $h_n: \Om(\Ga)\to \Om(\Ga_n)$ combined with the equivariant homeomorphism of the limit sets $\La(\Ga)\to \La(\Ga_n)$ yield an equivariant homeomorphism $\geo X\to \geo X$, see \cite{Tukia1985, Izeki}.   
Such an extension does not exist in higher rank since, in general, there is no equivariant homeomorphism of thickened limit sets $\Th(\Lat(\Ga))\to \Th(\Lat(\Ga_n))$. This can be already seen for  group actions on products of hyperbolic planes. 

(ii) An analogue of Theorem \ref{thm:Omstability} holds when we replace the group actions on flag-manifolds with actions on Finsler compactifications of the symmetric space and replace flag-manifold thickenings $\Th(\Lat)$  with Finsler thickenings $\ThF(\Lat)\subset \DF X$. Proving this requires extending Ehresmann--Thurston holonomy theorem to the category of smooth manifolds with corners and we will not pursue it here. 
\end{rem}

\subsection{Schottky actions}
\label{sec:schottky actions}

In this section we apply our local-to-global result for straight sequences 
(Theorem~\ref{thm:locstrimplcoastrseq})
to construct Morse actions of free groups, 
generalizing and sharpening\footnote{In the sense that we obtain free subgroups which are not only embedded, but also asymptotically embedded in $G$.} 
Tits's ping-pong construction.

\medskip 
We consider two oriented $\tau_{mod}$-regular geodesic lines 
$a,b$ in $X$. 
Let $\tau_{\pm a},\tau_{\pm b}\in\Flag(\tau_{mod})$ 
denote the simplices 
which they are $\tau$-asymptotic to, 
and let $\theta_{\pm a},\theta_{\pm b}\in\si_{mod}$ 
denote the types of their forward/backward ideal endpoints in $\geo X$.
(Note that $\theta_{-a}=\iota(\theta_a)$
and $\theta_{-b}=\iota(\theta_b)$.) Let $\Theta$ be a compact convex subset of $\ost(\tau_{mod})\subset \si_{mod}$, which is invariant under $\iota$.

\begin{dfn}[Generic pair of geodesics]
We call the pair of geodesics $(a,b)$ {\em generic} 
if the four simplices $\tau_{\pm a},\tau_{\pm b}$ 
are pairwise opposite. 
\end{dfn}
Let $\al,\beta\in G$ be axial isometries 
with axes $a$ and $b$ respectively 
and translating in the positive direction along these geodesics. 
Then $\tau_{\pm a}$ and $\tau_{\pm b}$ 
are the attractive/repulsive fixed points of $\al$ and $\beta$ 
on $\Flag(\tau_{mod})$. 

For every pair of numbers $m,n\in\N$ 
we consider the representation of the free group in two generators
\begin{equation*}
\rho_{m,n}: F_2=\<A,B\> \to G
\end{equation*}
sending the generator $A$ to $\al^m$ and $B$ to $\beta^n$. 
We regard it as an isometric action 
$\rho_{m,n}:F_2\acts X$.

\begin{definition}[Schottky subgroup]
A {\em $\tau_{mod}$-Schottky subgroup} of $G$ is a free $\tau_{mod}$-asymp\-to\-ti\-cally embedded subgroup of $G$. 
\end{definition}

If $G$ has rank one, this definition amounts to the requirement that $\Ga$ is convex cocompact and free. Equivalently, this is a discrete finitely generated subgroup of $G$ which contains  no nontrivial elliptic and parabolic elements and has totally disconnected limit set (see see \cite{Kapovich2007}). We note that this definition essentially agrees with the standard definition of Schottky groups in rank 1 Lie groups, provided one allows fundamental domains at infinity for such groups to be bounded by pairwise disjoint compact submanifolds which need not be topological spheres, see \cite{Kapovich2007} for the detailed discussion.

\begin{thm}[Morse Schottky actions]
\label{thm:mschott}
If the pair of geodesics $(a,b)$ is generic 
and if $\theta_{\pm a},\theta_{\pm b}\in\interior(\Theta)$, 
then the action $\rho_{m,n}$ is $\Theta$-Morse 
for sufficiently large $m,n$. Thus, such $\rho_{m,n}$ is injective and its image is a $\tau_{mod}$-Schottky subgroup 
of $G$. 
\end{thm}
\begin{rem}
In particular,
these actions are faithful and undistorted, 
compare Remark~\ref{rem:morseacundist}. 
\end{rem}
\proof
Let $S=\{A^{\pm1},B^{\pm1}\}$ be the standard generating set. 
We consider the sequences $(\ga_k)$ in $F_2$
with the property that $\ga_k^{-1}\ga_{k+1}\in S$ and 
$\ga_{k+1}\neq\ga_{k-1}$ for all $k$.  
They correspond to the geodesic segments in the Cayley tree of $F_2$ 
associated to $S$ which connect vertices. 

Let $x\in X$ be a base point. 
In view of Lemma~\ref{lem:morsefromhyp} 
we must show that the corresponding sequences $(\ga_kx)$ 
in the orbit $F_2\cdot x$
are uniformly $\Theta$-Morse. 
(Meaning e.g.\ that the maps $\R\to X$ 
sending the intervals $[k,k+1)$ to the points $\ga_kx$ are 
uniform $\Theta$-Morse quasigeodesics.)
As in the proof of Theorem~\ref{thm:locglobmqg}
we will obtain this by applying 
our local to global result for straight spaced sequences 
(Theorem~\ref{thm:locstrimplcoastrseq}) 
to the associated midpoint sequences. 
Note that the sequences $(\ga_kx)$ themselves cannot expected to be straight. 

Taking into account the $\Ga$-action, 
the uniform straightness of all midpoint sequences 
depends on the geometry of a finite configuration in the orbit. 
It is a consequence of the following fact.
Consider the midpoints $y_{\pm m}$ of the segments $x \al^{\pm m}(x)$ 
and $z_{\pm n}$ of the segments 
$x \beta^{\pm n}(x)$. 
\begin{lem}
\label{lem:geomquad}
For sufficiently large $m,n$ 
the quadruple $\{y_{\pm m},z_{\pm n}\}$ is arbitrarily separated 
and $\Theta$-regular.
Moreover, 
for any of the four points, 
the segments connecting it to the other three points 
have arbitrarily small $\zeta$-angles with the segment connecting it to $x$. 
\end{lem}
\proof
The four points are arbitrarily separated from each other and from $x$ 
because the axes $a$ and $b$ diverge from each other 
due to our genericity assumption. 

By symmetry, it suffices to verify the rest of the assertion 
for the point $y_m$, 
i.e.\ we show that 
the segments $y_my_{-m}$ and $y_mz_n$ are $\Theta$-regular 
for large $m,n$
and that 
$\lim_{m\to\infty} \angle^\zeta_{y_m}(x, y_{-m})=0$ 
and 
$\lim_{n,m\to\infty} \angle^\zeta_{y_m}(x, z_n)=0$. 

The orbit points $\al^{\pm m}x$ and the midpoints $y_{\pm m}$
are contained in a tubular neighborhood of the axis $a$. 
Therefore, 
the segments $y_mx$ and $y_my_{-m}$ 
are $\Theta$-regular for large $m$ 
and $\angle_{y_m}(x, y_{-m})\to0$. 
This implies that also $\angle^{\zeta}_{y_m}(x, y_{-m})\to0$. 

To verify the assertion for $(y_m,z_n)$ 
we use that, due to genericity, 
the simplices $\tau_a$ and $\tau_b$ are opposite 
and we consider the parallel set $P=P(\tau_a,\tau_b)$.
Since the geodesics $a$ and $b$ are forward asymptotic to $P$, 
it follows that the points $x,y_m,z_n$ have 
uniformly bounded distance from $P$. 
We denote their projections to $P$ by $\bar x,\bar y_m,\bar z_n$. 

Let $\Theta''\subset\interior(\Theta)$ be an auxiliary Weyl convex subset
such that 
$\theta_{\pm a},\theta_{\pm b}\in\interior(\Theta'')$. 
We have that 
$\bar y_m\in V(\bar x,\st_{\Theta''}(\tau_a))$ for large $m$ 
because the points $y_m$ lie in a tubular neighborhood 
of the ray with initial point $\bar x$ and asymptotic to $a$. 
Similarly, 
$\bar z_n\in V(\bar x,\st_{\Theta''}(\tau_b))$ for large $n$. 
It follows that 
$\bar x\in V(\bar y_m,\st_{\Theta''}(\tau_b))$ and,
using the convexity of $\Theta$-cones (Proposition~\ref{prop:thconeconv}), 
that $\bar z_n\in V(\bar y_m,\st_{\Theta''}(\tau_b))$. 

The cone $V(y_m,\st_{\Theta''}(\tau_b))$
is uniformly Hausdorff close to the cone 
$V(\bar y_m,\st_{\Theta''}(\tau_b))$
because the Hausdorff distance of the cones 
is bounded by the distance $d(y_m,\bar y_m)$ of their tips. 
Hence there exist points $x',z'_n\in V(y_m,\st_{\Theta''}(\tau_b))$ 
uniformly close to $x,z_n$. 
Since $d(y_m,x'),d(y_m,z'_n)\to\infty$ as $m,n\to\infty$, 
it follows that 
the segments $y_mx$ and $y_mz_n$ are $\Theta$-regular for large $m,n$. 
Furthermore, 
since $\angle^\zeta_{y_m}(x',z'_n)=0$ 
and 
$\angle_{y_m}(x,x')\to0$
as well as
$\angle_{y_m}(z_n,z'_n)\to0$, 
it follows that $\angle^\zeta_{y_m}(x,z_n)\to0$. 
\qed

\medskip
{\em Proof of Theorem concluded.} 
The lemma implies that for any given $l,\eps$ 
the midpoint triples of the four point sequences $(\ga_kx)$
are $(\Theta,\eps)$-straight and $l$-spaced 
if $m,n$ are sufficiently large, 
compare the quadruple condition (Definition~\ref{dfn:quad}). 
This means that the midpoint sequences 
of all sequences $(\ga_kx)$ are $(\Theta,\eps)$-straight and $l$-spaced 
for large $m,n$. 
Theorem~\ref{thm:locstrimplcoastrseq} then implies 
that the sequences $(\ga_kx)$ are uniformly $\Theta$-Morse. 
\qed 

\begin{rem}
1. Generalizing the above argument to free groups with finitely many generators, 
one can construct Morse Schottky subgroups 
for which the set $\theta(\La)\subset\si_{mod}$ of types of limit points 
is arbitrarily Hausdorff close to a given $\iota$-invariant 
Weyl convex subset $\Theta$. 
This provides an alternative approach 
to the second main theorem in \cite{Benoist}
using coarse geometric arguments. 

2. In \cite{DKL} Theorem \ref{thm:mschott} 
was generalized (by arguments similar to the its proof) to free products of Morse subgroups of $G$. 
\end{rem}

\subsection{Algorithmic recognition of Morse actions}
\label{sec:algrec}

In this section, we describe an algorithm which has an isometric action $\rho: \Ga\acts X$ and a point $x\in X$ 
as its input and terminates if and only if the action $\rho$ is Morse (otherwise, the algorithm runs forever). 

We begin by describing briefly the {\em Riley's algorithm} (see \cite{Riley}) accomplishing a similar task, namely, detecting 
geometrically finite actions on $X=\H^3$. 
Suppose that we are given a finite (symmetric) set of generators $g_1=1,\ldots, g_m$ 
of a subgroup  $\Ga\subset PO(3,1)$ and a base-point $x\in X=\H^3$. The idea of the algorithm is to construct a finite sided Dirichlet fundamental domain $D$ for $\Ga$ (with the center at $x$): Every geometrically finite subgroup of 
$PO(3,1)$ admits such a domain. (The latter is false for geometrically finite subgroups of $PO(n,1)$, $n\ge 4$, but is, nevertheless, true for convex cocompact subgroups.) Given  a finite sided convex fundamental domain, one concludes that $\Ga$ is geometrically finite. 
Here is how the algorithm works: For each $k$ define the subset $S_k\subset \Ga$ represented by words of length $\le k$ in  the letters $g_1,\ldots, g_m$. For each $g\in S_k$ consider the half-space $Bis(x, g(x))\subset X$ bounded by the bisector of the segment $x g(x)$ and containing the point $x$. Then compute the intersection
$$
D_k=\bigcap_{g\in S_k} Bis(x, g(x)).
$$ 
Check if $D_k$ satisfies the conditions of the {\em Poincar\'e's Fundamental Domain theorem}. 
If it does, then $D=D_k$ is a finite sided fundamental domain of $\Ga$. If not, increase $k$ by $1$ and repeat the process. Clearly, this process terminates if and only if $\Ga$ is geometrically finite. 

One can enhance the algorithm in order to detect if a geometrically finite group is convex cocompact. Namely, after a Dirichlet domain $D$ is constructed, one checks for the following:

1. If the ideal boundary of a Dirichlet domain $D$  has isolated ideal points (they would correspond to 
rank two cusps which are not allowed in convex cocompact groups). 

2. If the ideal boundary of $D$  contains tangent circular arcs with points of tangency fixed by parabolic elements (coming from the ``ideal vertex cycles''). Such points correspond to rank 1 cusps, which again are not allowed  
 in convex cocompact groups. 
 
Checking 1 and 2 is a finite process; after its completion, one concludes that $\Ga$ is convex cocompact.  
 
\medskip 
We refer the reader to \cite{Gilman1, Gilman2, Gilman-Maskit, Kapovich2015} and \cite[sect. 1.8]{manicures} for more details concerning discreteness algorithms  for groups acting on hyperbolic planes and hyperbolic 3-spaces.

\medskip
We now consider group actions on general symmetric spaces. Let $\Ga$ be a hyperbolic group with a fixed finite (symmetric) generating set; we equip the group $\Ga$ with the word metric determined by this generating set. 

For each $n$, let ${\mathcal L}_n$ denote the set of maps $q: [0, 3n]\cap \Z \to \Gamma$ which are restrictions  of geodesics $\tilde{q}: \Z\to \Ga$, such that $q(0)=1\in \Ga$. In view of the geodesic automatic structure on $\Ga$ (see e.g. \cite[Theorem 3.4.5]{Epstein}), the set  ${\mathcal L}_n$ can be described via a finite state automaton. 

Suppose that $\rho: \Gamma\acts X$ is an isometric action on a symmetric space $X$; we fix a base-point $x\in X$ 
and the corresponding orbit map $f: \Ga\to \Ga x\subset X$. We also fix an $\iota$-invariant face $\tau_{mod}$ of the model spherical simplex $\si_{mod}$ of $X$. The algorithm that we are about to describe will detect that the action $\rho$ is $\tau_{mod}$-Morse. 

\begin{rem}
If the face $\tau_{mod}$ is not fixed in advance, we would run algorithms for each face $\tau_{mod}$ in parallel. 
\end{rem}

For the algorithm we will be using a special (countable) increasing family of Weyl convex compact subsets $\Theta=\Theta_i\subset \ost(\tau_{mod})\subset \si_{mod}$ which exhausts $\ost(\tau_{mod})$; in particular, 
every compact $\iota$-invariant convex subset of $\ost(\tau_{mod})\subset \si_{mod}$ is contained in some $\Theta_i$: 
\begin{equation}\label{Theta_i}
\Theta_i:=\{v\in \si: \min_{\alpha\in \Phi_{\tau_{mod}}} \alpha(v)\ge \frac{1}{i} \},
\end{equation}
where $\Phi_{\tau_{mod}}$ is the subset of the set of simple roots $\Phi$ (with respect to $\si_{mod}$) which vanish on the face $\tau_{mod}$. Clearly, the sets $\Theta_i$ satisfy the required properties. Furthermore, we consider only those $L$ and $D$ which are natural numbers. 

Next, consider the sequence
$$
(L_i, \Theta_{i}, D_i)= (i, \Theta_i, D_i), i\in \N. 
$$

In order to detect $\tau_{mod}$-Morse actions we will use the local characterization of Morse quasigeodesics given by 
Theorem  \ref{thm:locstrimplcoastrseq} and Proposition \ref{prop:morseimplquad}. Due to the discrete nature of quasigeodesics that we will be considering, it suffices to assume that the additive quasi-isometry constant $A$ is zero. 

Consider the functions 
$$
l(\Theta, \Theta', \delta), \eps(\Theta, \Theta', \delta)
$$
as in Theorem \ref{thm:locstrimplcoastrseq}. Using these functions,  for the sets 
$\Theta=\Theta_{i}, \Theta'=\Theta_{i+1}$ and the constant $\delta=1$ we define the numbers 
$$
l_i=l(\Theta, \Theta', \delta), \eps_i= \eps(\Theta, \Theta', \delta). 
$$

Next, for the numbers $L=L_i, D=D_i$ and the sets  $\Theta=\Theta_{i}, \Theta'=\Theta_{i+1}$, consider the numbers 
$$
s_i=s(L_i,0, \Theta_i, \Theta_{i+1}, D_i, \eps_{i+1}, l_{i+1})
$$ 
as in Proposition \ref{prop:morseimplquad}. According to this proposition, every 
$(L_i,0,\Theta_i,D_i)$-Morse quasigeodesic satisfies the $(\Theta_{i+1},\eps_{i+1},l_{i+1}, s)$-quadruple 
condition for all $s\ge s_i$. We note that, a priori, the sequence $s_i$ need not  be increasing. We set $S_1=s_1$ and 
define a monotonic sequence $S_i$ recursively by 
$$
S_{i+1}= \max(S_i, s_{i+1}). 
$$
Then every $(\Theta_i, D_i, L_i, 0)$-Morse quasigeodesic also satisfies the 
$(\Theta_{i+1}, \eps_{i+1}, l_{i+1}, S_{i+1})$-quadruple condition. 

We are now ready to describe the algorithm. For each $i\in \N$ we compute the numbers $l_i, \eps_i$ and, then, $S_i$, as above.  
We then consider finite discrete paths in $\Ga$, $q\in {\mathcal L}_{S_i}$, and the corresponding discrete paths 
in $X$, $p(t)=q(t)x$, $t\in [0, 3S_i]\cap \Z$. The number of paths $q$ (and, hence, $p$) for each $i$ is finite, bounded by the growth function of the group 
$\Ga$. 

For each discrete path $p$ we check the $(\Theta_{i}, \eps_{i}, l_{i}, S_{i})$-quadruple condition. If for some $i=i_*$, all paths $p$ satisfy this condition, 
the algorithm terminates: It follows from Theorem \ref{thm:locstrimplcoastrseq} that the map $f$ sends all 
normalized discrete biinfinite geodesics in $\Ga$ to Morse quasigeodesics in $X$. Hence, the action $\Ga\acts X$ is 
Morse in this case. Conversely, suppose that the action of $\Ga$ is $(\Theta, D, L, 0)$-Morse. 
Then $f$ sends all isomeric embeddings $\tilde{q}: \Z\to \Ga$ to $(\Theta, D, L, 0)$-Morse quasigeodesics $\tilde{p}$ in $X$. 
In view of the properties of the sequence 
$$
(L_i, \Theta_i, D_i),
$$
it follows that for some $i$, 
$$
(L, \Theta, D)\le (L_i, \Theta_i, D_i), 
$$
i.e., $L\le L_i, \Theta\subset \Theta_i, D\le D_i$; hence, all the biinfinite discrete paths $\tilde p$ are 
$(\Theta_i, D_i, L_i, 0)$-Morse quasigeodesic. 
By the definition of the numbers $l_i, \eps_i, S_i$, it then follows that all the discrete paths $p=f\circ q, q\in {\mathcal L}_{S_i}$ satisfy the  
$(\Theta_{i+1}, \eps_{i+1}, l_{i+1}, S_{i+1})$-quadruple condition. Thus, the algorithm will terminate at the step $i+1$ in this case. 

Therefore, the algorithm terminates if and only if the action is Morse (for some parameters). If the action is not Morse, the algorithm will run forever. \qed

\begin{rem}
Applied to a rank one symmetric space $X$ and a hyperbolic group $\Ga$ without a nontrivial normal finite 
subgroup, the above algorithm verifies if the given representation 
$\rho: \Ga\to \Isom(X)$ is faithful with convex-cocompact image. 
We could not find this result in the existing literature; cf. however \cite{Gilman-Keen}. 
\end{rem}

\section{Appendix: Further properties of Morse quasigeodesics}

This is the only part of the paper not contained in \cite{morse}. Here we collect various properties of 
Morse quasigeodesics that we found to be useful elsewhere in our work.


\subsection{Finsler geometry of symmetric spaces}\label{sec:finsler}

In \cite{bordif}, see also \cite{anolec},
we considered a certain class of $G$-invariant ``polyhedral'' Finsler metrics on $X$.
Their geometric and asymptotic properties turned out to be well adapted 
to the study of geometric and dynamical properties of regular subgroups.
They provide a Finsler geodesic {\em combing} of $X$ which is, in many ways, more suitable for analyzing the asymptotic 
geometry of $X$ than the geodesic combing given by the standard Riemannian metric on $X$.
These Finsler metrics also play a basic role in the present paper.
We briefly recall their definition and some basic properties, 
and refer to \cite[\S 5.1]{bordif} for more details.

Let $\bar\theta\in\inte(\taumod)$ be a type spanning the face type $\taumod$. 
The {\em $\bar\theta$-Finsler distance} $d^{\bar\theta}$ on $X$ is the $G$-invariant pseudo-metric defined by
\begin{equation*}
d^{\bar\theta}(x,y) := \max_{\theta(\xi)=\bar\theta} \bigl( b_{\xi}(x)-b_{\xi}(y) \bigr) 
\end{equation*}
for $x,y\in X$, 
where the maximum is taken over all ideal points $\xi\in\geo X$ with type $\theta(\xi)=\bar\theta$.
It is positive, i.e. a (non-symmetric) metric, 
if and only if the radius of $\simod$ with respect to $\bar\theta$ is $<\pihalf$.
This is in turn equivalent to $\bar\theta$ not being contained in a factor of a nontrivial spherical join decomposition of $\simod$,
and is always satisfied e.g. if $X$ is irreducible.

If $d^{\bar\theta}$ is positive,
it is equivalent to the Riemannian metric. 
In general, if it is only a pseudo-metric,
it is still equivalent to the Riemannian metric $d$ on uniformly regular pairs of points.
More precisely,
if the pair of points $x,y$ is $\Theta$-regular, then 
$$
L^{-1} d(x,y) \le d^{\bar\theta}(x,y)\le L d(x,y)
$$ 
with a constant $L=L(\Theta)\ge  1$.

Regarding symmetry of the Finsler distance, one has the identity
\begin{equation*}
d^{\iota\bar\theta}(y,x)  = d^{\bar\theta}(x,y)
\end{equation*}
and hence $d^{\bar\theta}$ is symmetric if and only if 
$\iota\bar\theta=\bar\theta$.
We refer to $d^{\bar\theta}$ as a Finsler metric {\em of type $\taumod$}.

The $d^{\bar\theta}$-balls in $X$ are convex but not strictly convex.
(Their intersections with flats through their centers are polyhedra.)
Accordingly,
$d^{\bar\theta}$-geodesics connecting two given points $x,y$ are not unique. 
To simplify notation,
$xy$ will stand for {\em some} $d^{\bar\theta}$-geodesic connecting $x$ and $y$.
The union of all $d^{\bar\theta}$-geodesic $xy$ equals the $\taumod$-diamond $\diamot(x,y)$,
that is, a point lies on a $d^{\bar\theta}$-geodesic $xy$ if and only if it is contained in $\diamot(x,y)$, see \cite{anolec}.
Finsler geometry thus provides an alternative description of diamonds.  
Note that with this description, 
the diamond $\diamot(x,y)$ is also defined when the segment $xy$ is not $\taumod$-regular. 
Such a {\em degenerate} $\taumod$-diamond is contained in a smaller totally-geodesic subspace,
namely in the intersection of all $\taumod$-parallel sets containing the points $x,y$.
The description of geodesics and diamonds 
also implies that the unparameterized $d^{\bar\theta}$-geodesics depend only on the face type $\taumod$,
and not on $\bar\theta$.
We will refer to $d^{\bar\theta}$-geodesics as {\em $\taumod$-Finsler geodesics}.
Note that Riemannian geodesics are Finsler geodesics.

We will call a $\Theta$-regular $\taumod$-Finsler geodesic a {\em $\Theta$-Finsler geodesic}.
If $xy$ is a $\Theta$-regular (Riemannian) segment,
then the union of $\Theta$-Finsler geodesics $xy$ equals the $\Theta$-diamond $\diamoTh(x,y)$. 

Every $\taumod$-Finsler ray in $X$ is contained in a $\taumod$-Weyl cone,
and we will use the notation $x\tau$ for a $\taumod$-Finsler ray contained $V(x,\st(\tau))$.
Similarly, 
every $\taumod$-Finsler line is contained in a $\taumod$-parallel set,
and we denote by $\tau_-\tau_+$ an oriented $\taumod$-Finsler line 
forward/backward asymptotic to two antipodal simplices $\tau_{\pm}\in\Flagt$ 
and contained in $P(\tau_-,\tau_+)$.

Examples of $\Theta$-regular Finsler geodesics can be obtained as follows. Let $(x_i)$ be a (finite or infinite) 
sequence contained in a parallel set $P(\tau_-,\tau_+)$ such that each Riemannian segment $x_i x_{i+1}$ is 
$\tau_+$-longitudinal and $\Theta'$-regular. Then the concatenation of these geodesic segments is 

Conversely, every $\Theta$-regular Finsler geodesic $c: I\to X$ can be {\em approximated} by a piecewise-Riemannian Finsler geodesic $c'$: 
Pick a number $s>0$ and consider a maximal $s$-separated subset $J\subset I$. Then take $c'$ to be the concatenation of Riemannian geodesic 
segments $c(i)c(j)$ for consecutive pairs $i, j\in J$. In view of this approximation procedure, the String of Diamonds Theorem (Theorem \ref{thm:string}) 
holds if instead of Riemannian geodesic segments $x_i x_{i+1}$ we allow $\Theta$-regular Finsler segments.

\subsection{Stability of diamonds}

Diamonds can be regarded as Finsler-geometric replacements
of {\em geodesic segments} in nonpositively curved symmetric spaces of higher rank.

Riemannian geodesic segments in Hadamard manifolds (and, more generally, $CAT(0)$ metric spaces) 
depend {\em uniformly continuously} on their tips: By convexity of the distance function we have,  
$$
d_{Haus}(xy, x'y')\le \max( d(x, x'), d(y, y')). 
$$ 
In \cite[Prop. 3.70]{mlem} we proved that diamonds $\diamot$  depend  {\em continuously} on their tips. 

Below we establish uniform control on how much sufficiently large $\Theta$-diamonds vary with their tips. 

\begin{lem}
\label{lem:ucldiamonds}
For $d'>d>0$ there exists $C=C(\Theta,\Theta',d,d')$ such that the following holds:

If a segment $x_-x_+\subset X$ is $\Theta$-regular with length $\geq C$
and $y_{\pm}\in B(x_{\pm},d)$,
then the segment $y_-y_+$ is $\Theta'$-regular and 
$\diamoTh(x_-, x_+)\subset N_{d'}(\diamoThp(y_-, y_+))$.
\end{lem}
\proof 
The $\Theta'$-regularity of $y_-y_+$ for sufficiently large $C$ follows from the $\De$-triangle inequality.

Suppose that there exists no constant $C$ 
for which also the second assertion holds. 
Then there are sequences of points $x_n^{\pm}$ with $d(x_n^-,x_n^+)\to+\infty$,
$y_n^{\pm}$ with $d(x_n^{\pm},y_n^{\pm})\leq d$, 
$x_n\in\diamoTh(x_n^-,x_n^+)$ and $y_n\in\diamoThp(y_n^-,y_n^+)$
with $d(x_n,\diamoThp(y_n^-,y_n^+))=d(x_n,y_n)=d'$.
We may assume convergence $x_n\to x_{\infty}$ and $y_n\to y_{\infty}$ in $X$.

After extraction, at least one of the sequences $(x_n^{\pm})$ diverges. 
There are two cases to consider.

Suppose first that both sequences $(x_n^{\pm})$ diverge.
Then they are uniformly $\taumod$-regular and, 
after extraction,
we have $\taumod$-flag convergence $x_n^{\pm},y_n^{\pm}\to\tau_{\pm}\in\Flagt$.
The limit simplices $\tau_{\pm}$ are antipodal (because $x_n\to x_{\infty}$).
We observe that 
$$
d(x_n,\D\diamoThp(x_n^-,x_n^+)),d(y_n,\D\diamoThp(y_n^-,y_n^+))\to+\infty.$$ 
It follows that the sequences of diamonds $\diamoThp(x_n^-,x_n^+)$ and $\diamoThp(y_n^-,y_n^+)$ 
both Hausdorff converge to the $\taumod$-parallel set $P=P(\tau_-,\tau_+)$. 
It holds that $x_{\infty}\in P$ because $x_n\in \diamoTh(x_n^-,x_n^+)$.
On the other hand, $d(x_{\infty},P)=d'$ because $d(x_n,\diamoThp(y_n^-,y_n^+))=d'$,
a contradiction. 

Second, suppose that only one of the sequences $(x_n^{\pm})$ diverges, say,
after extraction, 
$x_n^-\to x_{\infty}^-$ and $y_n^-\to y_{\infty}^-$ in $X$ to limit points with $d(x_{\infty}^-,y_{\infty}^-)\leq d$,
and $x_n^+\to\tau_+\in\Flagt$.
Now the distance of $x_n$ from the boundary of the $\Theta'$-Weyl cone with tip $x_n^+$ and containing $x_n$ 
goes to infinity and it follows that 
$\diamoThp(x_n^-,x_n^+)\to V(x_{\infty}^-,\st_{\Theta'}(\tau_+))$ and, similarly, 
$\diamoThp(y_n^-,y_n^+)\to V(y_{\infty}^-,\st_{\Theta'}(\tau_+))$.
The asymptotic limit Weyl cones have Hausdorff distance $d(x_{\infty}^-,y_{\infty}^-)$.
On the other hand, 
$x_{\infty}\in V(x_{\infty}^-,\st_{\Theta'}(\tau_+))$
and $d(x_{\infty},V(y_{\infty}^-,\st_{\Theta'}(\tau_+)))=d'$, again a contradiction. 

This shows that also (ii) holds for sufficiently large $C$.
\qed 

\medskip
We reformulate this result in terms of Finsler geodesics:
\begin{lem}
\label{lem:uclfgeos}
There exists $C=C(\Theta,\Theta',d,d')$ such that the following holds: 
If $x_-x_+$ is a $\Theta$-Finsler geodesic in $X$ with $d(x_-,x_+)\geq C$
and $y_{\pm}$ are points with $d(y_{\pm},x_{\pm})\leq d$,
then every point $x$ on $x_-x_+$ lies within distance $d'$ of a point $y$ on a $\Theta'$-Finsler geodesic $y_-y_+$. 
\end{lem}

Note that we do not claim here that one can take the same Finsler geodesic $y_-y_+$ for all points $x$ on $x_-x_+$.

\medskip 
We now apply this stabilty result to Morse quasigeodesics. One, somewhat annoying, feature of the definition of  
$\Theta$-Morse quasigeodesics $p: I\to X$ is that $p([t_1,t_2])$ is not required to be uniformly close to a $\Theta$-diamond spanned 
by $p(t_1), p(t_2)$.  (One reason is  because the segment $p(t_1)p(t_2)$ need not be $\Theta$-regular.)  
Nevertheless, Lemma \ref{lem:ucldiamonds} implies: 

\begin{lemma}\label{lem:Morse-close-to-diamond}
For every Morse datum $M=(\Theta,B,L,A)$ and $\Theta'> \Theta$, there exists $C=C(M,\Theta')$ and $D'$ 
such that whenever $d(x_1, x_2)\ge C$, the segment  $x_1x_2= p(t_1)p(t_2)$ is $\Theta'$-regular and 
$p([t_1,t_2])$ lies in the  $D'$-neighborhood of the $\Theta'$-diamond $\diamo_{\Theta'}(x_1, x_2)$.  
\end{lemma}

\subsection{Finsler approximation of Morse quasigeodesics}


The next theorem establishes that every (sufficiently long) Morse quasigeodesic is uniformly close to a Finsler geodesic with 
the same end-points. In this theorem, for convenience of the notation, 
we will be allowing Morse quasigeodesics $p$ to be defined on closed intervals $I$ in the extended real line; 
this is just a shorthand for a map $I'=I\cap \R\to X$ such that, as $I'\ni t\to \pm \infty$, $p(t)\to p(\pm \infty)\in \Flagt$.  
When we say that such maps $p, c$ are within distance $D'$ from each other, this simply means that their restrictions to 
$I'$ are within distance $\le D'$.

\begin{thm}
[Finsler approximation theorem]  \label{thm:Finsler-approximation} 
For every Morse datum $M=(\Theta, D, L, A)$, $\Theta'> \Theta$, and a positive number $S$, there exist  
$C=C(M, \Theta', S), D'=D'(M, \Theta', S)$ satisfying  the following.  
 
Let $p: I=[t_-, t_+]\to X\cup \Flagt$ be a $M$-Morse quasigeodesic between the points $x_\pm=p(t_\pm)\in X\cup \Flagt$ such that 
$d(x_-, x_+)\ge C$.  Then  there exists a $\Theta'$-Finsler geodesic $x_-x_+$ equipped with a monotonic parameterization 
$c: I \to x_-x_+$ such that:

(a) The maps $p, c: I\to X$ are within distance $\le D'$ from each other. 

(b)  $x_-x_+$ is an $S$-spaced piecewise-Riemannian geodesic, i.e.  the Riemannian 
length of each Riemannian segments of  $x_-x_+$ is $\ge S$. 
\end{thm}
\proof  We will prove this in the case when both $x_\pm$ are in $X$ since the proofs when one or both points 
$x_\pm$ are in $\Flagt$ are similar: One replaces diamonds with Weyl cones or parallel sets.

By the definition of an $M$-Morse quasigeodesic,  for all subintervals $[s_-, s_+]\subset [t_-, t_+]$,  
there exists a $\Theta$-diamond 
$$
\diamoTh(y'_-, y'_+)
$$ 
whose $D$-neighborhood contains $p([s_-, s_+])$, and for $y_\pm=p(s_\pm)$, we have 
$$
d(y_\pm,  y'_\pm)\le D. 
$$
Therefore, applying the first part of 
Lemma \ref{lem:ucldiamonds}, we conclude that the Riemannian segment $y_- y_+$  is $\Theta'$-regular provided that  
$d(y_-, y_+)\ge C_1=C_1(M,\Theta')$.  In view of the quasigeodesic property of $p$, the last inequality follows from the separation condition
$$
s_+ - s_- \ge s=s(M,\Theta'). 
$$
This, of course, also applies to $[s_-, s_+]= [t_-, t_+]$ and, hence, using 
the second part of Lemma  \ref{lem:ucldiamonds}, we obtain   
$$
p(I)\subset N_D\left( \diamoTh (x'_-, x'_+) \right) \subset N_{D+D_1} \left( \diamoThp (x_-, x_+) \right), 
$$
where $D_1=D_1(M,\Theta')$. We let 
$$
\bar{y}_\pm\in \diamo':= \diamoThp (x_-, x_+)= V(x_-, \st_{\Theta'}(\tau_+))\cap V(x_+, \st_{\Theta'}(\tau_-))
$$
denote the nearest-point projections of $y_\pm=p(s_\pm)$. 
As long as $s_+ - s_-\ge s'(M,\Theta')$, the Riemannian segments $\bar{y}_- \bar{y}_+$ are also 
$\Theta'$-regular and have length $\ge S$.  Furthermore, as in the proof of Proposition  \ref{prop:morseimplquad}, we can choose 
$s'$  such that each segment $\bar{y}_- \bar{y}_+$ is $\tau_+$-longitudinal. 

We assume, from now on, that  $t_+ - t_-\ge s''(M,\Theta')$,  which is achieved by assuming that
$$
L^{-1} (d(x_-, x_+) - A) \ge s'(M,\Theta'). 
$$
Take  a  maximal $s'$-separated  subset  $J\subset I$ containing $t_\pm$.  For each $j\in J$ define the point 
$$
z_j:= \overline{p(j)} \in \diamo'. 
$$
Then for  all consecutive $i, j\in J$, $s'\le |j-i|\le 2s'$ we have  
\begin{equation}\label{eq:double}
L^{-1} s' - (A+2D+2D_1)  \le d(z_i, z_j) \le 2Ls' + (A+2D+2D_1) .  
\end{equation}
We then let $c$ denote the concatenation of  Riemannian segments $z_{i} z_{j}$ for consecutive $i, j\in J$,  
where we use the affine parameterization of  $[i,j]\to z_{i} z_{j}$.  Thus, $c$ is a $\Theta'$-Finsler geodesic. 
We now take the smallest $s''\ge s'(M,\Theta')$ satisfying 
 $$
S\le  L^{-1} s'' - (A+2D+2D_1),
$$
the inequalities \eqref{eq:double} imply that $c$ satisfies both requirements of the approximation theorem with
$$
D'= 2Ls'' + (A+2D+2D_1) + (D+D_1) + (2Ls'' +A). \qed 
$$

\begin{rem}
In the case when the domain of $p$ is unbounded, one can prove a bit sharper result, namely, one can take $\Theta'=\Theta$. Compare \cite[sect. 6]{relmorse-1}. 
\end{rem}

\subsection{Altering Morse quasigeodesics}

\medskip 
Below we consider certain modifications of $M$-Morse quasigeodesics $p$ in $X$ represented as concatenations 
$p=p_-\star p_0 \star p_+$, where $x_\pm$ are the end-points of $p_0$, and $y_\pm, x_\pm$ are the end-points of $p_\pm$. 
(As in the previous section, we will be allowing $y_\pm$ to be in $X\cup \Flagt$.)   These modifications will have the form 
$p'= p'_- \star p'_0 \star p'_+$, where $p'_\pm$ and $p'_0$ are all Morse. 
We will see that, under certain assumptions, the entire  $p'$ is again Morse  (for suitable Morse datum $M'$). 

We begin by analyzing extensions of $p$ to biinfinite paths. 

\begin{lemma}
[Extension lemma]   \label{lem:extension}
Suppose that 
$$
p_\pm \subset V_\pm= V(x_\pm, \st(\tau_\pm)). 
$$
Whenever $y_\pm$ is in $X$,  we let  $c_\pm$ be $\Theta$-regular Finsler rays contained in $V_\pm$ and connecting $y_\pm$ to $\tau_\pm$. 
Then, for every $\Theta'> \Theta$, there exists a Morse datum $M'$ containing $\Theta'$  such that the concatenation
$$
\hat{p}= c_- \star p\star c_+
$$
is $M'$-Morse, provided that $d(x_\pm, y_\pm)\ge C=C(M, \Theta')$.
\end{lemma}
\proof  We fix an auxiliary subset $\Theta_1$ satisfying $\Theta< \Theta_1 < \Theta'$. 
We let $S=S(\Theta_1,\Theta',1), \eps=\eps(\Theta_1,\Theta',1)$ be  constants as in the string of diamonds theorem 
(Theorem \ref{thm:string}). 

According to Theorem \ref{thm:Finsler-approximation}, there exists a $\Theta'$-regular Finsler geodesic 
$$
\bar{c}= y_- \bar{x}_- \star \bar{x}_- \bar{x}_+ \star  \bar{x}_+ y_+
$$
within distance $D_1=D_1(M,\Theta',S)$ from the path $p$, such that $\bar{c}$ 
is the concatenation of   segments of length $\ge S$ and $d(x_\pm, \bar{x}_\pm)\le D_1$.  
We let $z_\pm y_\pm$ denote the  subsegments of $\bar{x}_\pm y_\pm$ containing $y_\pm$.

Since $d(x_\pm, \bar{x}_\pm)\le D_1$, for each $\eps>0$ and 
a sufficiently large $C_1=C_1(D_1,\Theta')$,  the inequality $d(x_\pm, y_\pm)\ge C_1$ implies  
$$
\zangle_{y_{\pm}}(x_\pm, \bar{x}_\pm)\le \eps.
$$
Therefore,
$$
\zangle_{y_\pm}(z_\pm, \tau_\pm )\ge \pi -\eps 
$$
and, hence, the piecewise-geodesic path
$$
\hat{c}=c_- \star \bar{c} \star c_+
$$
is $(\Theta_1,\eps)$-straight and $S$-spaced. Hence, by Theorem \ref{thm:string}, the concatenation $\hat{c}$ is $M'$-Morse, 
where $M_1=(\Theta', 1, L, A)$. Since the path $\hat{p}$ is within distance $D_1$ from $\hat{c}$, it is $M'$-Morse, where 
$M'= M_1 + D_1$. \qed  

\medskip 
The next lemma was proven in \cite[Thm. 4.11]{DKL} in the case when $p, p'$ are finite paths. 
The proof in the case of (bi)infinite paths is the same and we omit it.  

\begin{lemma}[Replacement lemma] 
\label{lem:replacement}
Suppose that $p'= p'_- \star p'_0 \star p'_+$ is a concatenation of $M$-Morse quasigeodesics in $X$, such that the end-points of $p_\pm, p'_\pm$ and $p_0, p'_0$ are the same.  Then for every $\Theta'> \Theta$ there exists a Morse datum $M'$ containing $\Theta'$ such that the path $p'$ is $M'$-Morse.  
\end{lemma}

In the following lemmata we will modify the path $p$ by altering $p_\pm$ and keeping $p_0$ unchanged or moving it by a small amount 
(``wiggling the head and the tail of $p$''). 

\begin{lemma}
[Wiggle lemma, I]\label{lem:wiggle1}
 Suppose that the paths $p_\pm, p'_\pm$ are both infinite. We let $p'_\pm$ be $M$-Morse quasigeodesics with finite terminal points $x_\pm$ and set 
$p':= p'_-\star p_0 \star p'_+$.  Given $\Theta'> \Theta$ there exists $\eps=\eps(M, \Theta')>0$ and a Morse datum $M'$ containing $\Theta'$  
such that if 
$$
\mu:= \max( \zangle_{x_\pm}(p'_\pm(\pm \infty),  p_\pm(\pm\infty)) ) < \eps,
$$
then $p'$ is $M'$-Morse. 
\end{lemma}
\proof  We fix an auxiliary compact Weyl-convex subset $\Theta_1\subset \ost(\taumod)$ such that $\Theta< \Theta_1  < \Theta'$.  
Set $\tau_\pm= p_\pm(\pm\infty)$, $\tau'_\pm= p'_\pm(\pm\infty)$. 

According to Lemma \ref{lem:replacement}, there exists a Morse datum $M_1$ containing $\Theta_1$ such that for any 
$\Theta_1$-regular Finsler geodesic rays $c_\pm:= x_\pm \tau_\pm$, the concatenation $c_-\star p_0 \star c_+$ is $M_1$-Morse.

Let $M_2 > M_1+1$ be a Morse datum containing $\Theta'$ and 
let $S>0$ be such that if a path $q$ in $X$ is $S$-locally $M_1+1$-Morse then $q$ is $M_2$-Morse (see Theorem  \ref{thm:locglobmqg}).   
Let $\eps$ be such that for $x\in X, \tau, \tau'\in \Flagt$, if $\zangle_x(\tau, \tau')<\eps$  then each $\Theta_1$-regular Finsler segment of length $\le S$ in $V(x, \st(\tau'))$ is 
within unit distance from a $\Theta_1$-regular Finsler segment of length $\le S$ in $V(x, \st(\tau))$.  We assume now that $\mu<\eps$.

Since $p'_\pm$ are $M$-Morse rays, they are within distance $D_1=D_1(M,\Theta_1)$ from $\Theta_1$-regular Finsler   rays $c'_\pm= x_\pm \tau'_\pm$ 
connecting  $x_\pm$ and $\tau'_\pm$.  Define a new path $c':= c'_-\star p_0 \star c'_+$.

By our choice of $\eps$, the  $\Theta_1$-regular Finsler subsegment $s'_\pm = x_\pm y'_\pm$ of $c'_\pm$  of length $S$ is within unit distance from a 
 $\Theta_1$-regular Finsler subsegment $s_\pm = x_\pm y_\pm$ of $c_\pm$   of length $S$, where $c_\pm= x_\pm \tau_\pm$ 
is a $\Theta_1$-Finsler geodesic connecting $x_\pm$ to $\tau_\pm$.

The concatenation 
$$
s_- \star p_0 \star s_+
$$
is $M_1$-Morse,  and, since $c'_\pm$ are $\Theta_1$-Finsler geodesic, the path  $c'$ is $S$-locally $M_1+1$-Morse.  
By our choice of $S$, the path $c'$ is $M_2$-Morse.  
Since $c'$ is within distance $D_1$ from $p'$, the path  $p'$ is $M_2+D_1$-Morse.  Lastly, we set $M':= M_2 +D_1$. \qed

\medskip 

We generalize this lemma by allowing finite Morse quasigeodesics. 
 We continue with the setting of Lemma \ref{lem:wiggle1}; we now allow paths $p_\pm$ and $p'_\pm$ to be finite, connecting $y_\pm, x_\pm$ and $y'_\pm, x_\pm$ respectively. (Some of $y_\pm, y'_\pm$ might be in $\Flagt$.) 
However, we will assume that the distances $d(x_\pm, y_\pm), d(x'_\pm, y_\pm)$ are sufficiently large, $\ge C$.

\begin{lemma}
[Wiggle lemma, II]\label{lem:wiggle2}
 Given $\Theta'> \Theta$ there exist $C\ge 0$, $\eps>0$   and a Morse datum $M'$ containing $\Theta'$  such that   if  
$$
\mu:= \max( \zangle_{x_\pm}(y'_\pm,  y_\pm) ) < \eps,
$$
and
$$
\nu:= \min( d(x_\pm, y_\pm), d(x_\pm, y'_\pm)   ) \ge C 
$$
then $p'$ is $M'$-Morse.  
\end{lemma}
\proof  Pick  an auxiliary compact Weyl-convex subset $\Theta_2$, $\Theta< \Theta_2 < \Theta'$.

We define biinfinite geodesic extensions $\hat{p}, \hat{p}'$ as in Lemma \ref{lem:extension}, by extending (if necessary) the paths $p_\pm, {p}'_\pm$  via $\Theta$-Finsler  geodesics 
$y_\pm \tau_\pm$ and $y'_\pm \tau'_\pm$.  According to Lemma \ref{lem:extension},  there exists $C> 0$and a Morse datum $M_2$ (containing $\Theta_2$), 
both  depending on $M$ and $\Theta_2$,  such that  the path $\hat{p}$ is $M_2$-Morse. The  same lemma applied to 
the paths $\hat{p}'_\pm$ implies that they are also $M_2$-Morse. 

By the construction,  
$$
\mu:= \zangle_{x_\pm}(y'_\pm,  y_\pm) = \zangle_{x_\pm}(\tau'_\pm,  \tau_\pm).  
$$
 Now, claim follows from Lemma \ref{lem:wiggle1}. \qed 

\medskip 
Lastly, we prove a general Wiggle Lemma where we allow to perturb the entire path $p$. We consider concatenations 
$$
p= p_- \star p_0 \star p_+, \quad p'= p'_- \star p'_0 \star p'_+
$$
of $M$-Morse quasigeodesics, where we assume that $p_0, p'_0$ are within distance $D_0$ from each other.  The paths $p_\pm$ connect $y_\pm, x_\pm$ and 
$p'_\pm$ connect $y'_\pm, x'_\pm$.

\begin{lemma}
[Wiggle lemma, III]\label{lem:wiggle3}
 Given $\Theta'> \Theta$ there exist $C\ge 0$, $\eps>0$  and a Morse datum $M'$ containing $\Theta'$  such that   if  
$$
\mu:= \max( \zangle_{x_\pm}(y'_\pm,  y_\pm) ) < \eps,
$$
and
$$
\nu:= \min( d(x_\pm, y_\pm), d(x'_\pm, y'_\pm)   ) \ge C 
$$
then $p'$ is $M'$-Morse.  
\end{lemma}
\proof As before, we fix an auxiliary compact Weyl-convex subset $\Theta_3$, $\Theta< \Theta_3 < \Theta'$.  Then $p'_\pm$ are within distance 
$D_3=D_3(M, \Theta_3)$ from  $\Theta_3$-regular Finsler geodesics $c_\pm:= y'_\pm x_\pm$.  We  apply Lemma  \ref{lem:wiggle2} to the pair of paths 
$$
p, p'':= c_- \star p_0 \star c_+. 
$$
It follows that $p''$ is $M_3$-Morse for some Morse datum $M_3$ containing $\Theta'$ provided that $\mu\le \eps= \eps(M, \Theta_3, \Theta')$ 
and $\nu \ge C=C(M, \Theta_3, \Theta')$. Since 
the paths $p''$ and $p'$ are wihin distance $D':= \max(D_0,D_3)$ from each other, the path $p'$ is $M':=M_3+D'$-Morse. \qed

Addresses:

\noindent M.K.: Department of Mathematics, \\
University of California, Davis\\
CA 95616, USA\\
email: kapovich@math.ucdavis.edu

\noindent B.L.: Mathematisches Institut\\
Universit\"at M\"unchen \\manicures
Theresienstr. 39\\ 
D-80333, M\"unchen, Germany\\ 
email: b.l@lmu.de

\noindent J.P.: Departament de Matem\`atiques,\\
 Universitat Aut\`onoma de Barcelona,\\ 
 08193 Bellaterra, Spain\\
email: porti@mat.uab.cat

\end{document}